\documentstyle[12pt]{article}
\input amssymb.sty
\newcommand{\bbbc}{{\Bbb C}}
\newcommand{\mn}{\medskip\noindent}
\newcommand{\bn}{\bigskip\noindent}
\newcommand{\sn}{\smallskip\noindent}
\newcommand{\X}{{\cal{X}}}
\newcommand{\A}{{\cal{A}}}
\newcommand{\T}{{\cal{T}}}
\newcommand{\id}{{\rm{id}}}

\newcommand{\Y}{{\cal{Y}}}
\newcommand{\Z}{{\cal{Z}}}
\newcommand{\cO}{{\cal{O}}}
\newcommand{\kd}{{\rm d}}
\newcommand{\r}{{\rm r}}
\newcommand{\br}{{\bf r}}
\newcommand{\uu}{{\bf u}}
\newcommand{\Gammalinv}{_{\rm inv}\Gamma}
\newcommand{\llabel}[1]{\label{#1}}
\begin{document}

\title{On the Construction of Covariant Differential Calculi on 
Quantum Homogeneous Spaces}
\date{}
\maketitle
\vspace*{-8ex}
\centerline{Konrad Schm\"udgen}
\centerline{Fakult\"at f\"ur Mathematik und Informatik}
\centerline{Universit\"at Leipzig, Augustusplatz 10, 04109 Leipzig, Germany}
\centerline{E-mail: schmuedg@mathematik.uni-leipzig.de}

\begin{abstract} Let $\A$ be a coquasitriangular Hopf algebra and 
$\X$ the subalgebra of $\A$ generated by a row of a matrix 
corepresentation $\uu$ or by a row of $\uu$ and a row of the 
contragredient representation $\uu^c$. In the paper 
left-covariant first order differential 
calculi on the quantum group $\A$ are constructed and the corresponding 
induced calculi on the left quantum space $\X$ are described. 
The main tool for these constructions are the L-functionals 
associated with $\uu$. The results are applied to the quantum 
homogeneous space 
$GL_q(N) / GL_q(N{-}1)$. 
\end{abstract}
\renewcommand{\baselinestretch}{1.0}
\setcounter{section}{-1}
\section{Introduction}
Based on the pioneering work of S.L. Woronowicz [W2], a beautiful 
theory of bicovariant differential theory on 
quantum groups has been developed till now. A thorough treatment of 
this theory can be found in Chapter 14 of the monograph [KS]. The 
theory of covariant differential calculi
on quantum spaces, in contrast, is still at the very beginning and 
neither general methods 
for the construction of such calculi nor remarkable 
general results are known. Covariant differential calculi have been constructed and 
studied so far only on a few simple quantum spaces 
([PW], [WZ], [P1], [P2], [SS1], [AS], [CHZ], [We] ). 

In this paper we are concerned with the construction of first order 
differential calculi (abbreviated, FODC) on subalgebras of a coquasitriangular 
Hopf algebra $\A$ which are generated by a row of a 
fixed corepresentation $\uu$ or by a row of $\uu$ 
and a row of the contragredient corepresentation $\uu^c$ of $\A$. Such a subalgebra is a left quantum space of $\A$ with left coaction 
given by the restriction of the comultiplication. Our method of 
construction is easy to explain: The FODC on the quantum spaces 
are induced from appropriate {\it left-covariant} differential 
calculi on the quantum group $\A$. The main technical tool for the construction 
of the left-covariant calculi on $\A$ are the L-functionals associated with 
the corepresentation $\uu$. We always try to be as simple and 
close to the classical situation as possible. Our approach has two important 
advantages: First, because of the close relationship between the 
calculi on the quantum space and on the quantum group the theory of 
L-functionals and other Hopf algebra techniques can be applied to the 
study of the 
calculi on the quantum space. Secondly, the simplicity of the 
constructed left-covariant calculi, in contrast to the usual 
bicovariant calculi, might be useful for doing explicit 
computations. Our guiding example are the quantum 
spheres associated with the quantum group $GL_q(N)$ (see [SV], [NYM] or 
[KS], 11.6). For these quantum spheres a classification of covariant 
differential calculi has been recently given by M. Welk [We]. As 
an application of our method we describe 
some of the main calculi occuring there as induced 
from left-covariant calculi on $GL_q(N)$. Strictly speaking, we derive 
the left-covariant counter-parts of these calculi, 
because in [We] right quantum spheres and 
right-covariant calculi are investigated.

This paper is organized as follows. Section 1 contains 
some preliminaries and collects some notation. 
In Sections 2 and 3 first order calculi on the left 
quantum spaces generated a a single row of $\uu$ and $\uu^c$, respectively, are 
investigated. Section 4 deals with the quantum space generated by a row of $\uu$ 
and a row of $\uu^c$. Four families of covariant FODC are constructed and 
the commutation rules between generators and their differentials are 
explicitely described. The application of the results to the fundamental 
corepresentations of the quantum groups $GL_q(N)$ and $SL_q(N)$ 
are discussed in Section 5. In Section 6 another interesting 
FODC on the quantum sphere is obtained 
from a particular {\it bicovariant} (!) calculus on $GL_q(N)$. The left-covariant 
differential calculi on the quantum groups have been so far only auxilary 
tools for the study of the induced FODC on the quantum spaces. In Section 7 
the same idea is used in order to construct "reasonable" 
left-covariant FODC on the quantum groups 
$GL_q(N), SL_q(N), O_q(N)$ and $Sp_q(N)$ which are in many aspects 
close to the ordinary differential calculus on the corresponding 
Lie groups. In particular, the dimensions of these calculi 
coincide with the classical group dimensions. 

I would like to thank M. Welk for useful 
discussions on the subject of the paper.

\section*{2\ Preliminaries}
Throughout this paper $\A$ is a coquasitriangular 
complex Hopf algebra and $\br$ denotes a fixed universal r-form of $\A$ 
(see, for instance, [LT] or [KS], Section 10.1, for these notions). 
The comultiplication, the counit and the antipode of $\A$ are 
denoted by $\Delta$  $\varepsilon$ and $S$, respectively. 
We shall use the
Sweedler notation $\Delta(a)=a_{(1)}\otimes a_{(2)}$ 
for the comultiplication of $\A$. Let us recall that a 
Hopf algebra $\A$ is called {\it coquasitriangular} if it is equipped with 
a linear functional $\br$ on $\A \otimes \A$ which 
is invertible with respect to the convolution multiplication and satisfies the 
following conditions for arbitrary elements $a,b,c \in \A$:
\begin{eqnarray}\llabel{Kr0}
\br( ab \otimes c)= \br( a \otimes c_{(1)})\br( b\otimes c_{(2)}),~~
\br( a\otimes bc)= \br( a_{(1)}\otimes c)\br( a_{(2)}\otimes b),\\
\llabel{Kr1}
\br( a_{(1)}\otimes b_{(1)}) a_{(2)} b_{(2)}=\br( a_{(2)}\otimes
b_{(2)}) b_{(1)} a_{(1)}.\qquad\qquad
\end{eqnarray}
Such a linear form $\br$ is called a 
{\it universal $r$-form} of the Hopf algebra $\A$. The convolution 
inverse of $\br$ is denoted by $\bar{\br}$. We shall write $\br (a,b) :=
\br (a \otimes b), a,b \in \A$.

Further, $\uu =(u^i_j)_{i,j=1,\ldots,n}$  denotes a fixed $n$-dimensional 
{\it matrix corepresentation} of $\A$, that is, $\uu$ is an 
$n \times n$-matrix of elements $u^i_j$ of $\A$ 
such that 
$$\Delta (u^i_j) = \sum_{k=1}^n u^i_k \otimes u^k_j ~ \mbox { and }
\varepsilon (u^i_j) = \delta_{ij}~ \mbox{ for}~ i,j=1,\cdots, n.
$$ 
We define the L-functionals ${l^{\pm}}^i_j$ and the R-matrix 
$\hat{R}$ associated with the 
corepresentation $\uu$ by
\begin{eqnarray*}
{l^+}^i_j(\cdot)=\br(\cdot\otimes u^i_j),~~
{l^-}^i_j(\cdot)=\bar{\br}(u^i_j\otimes\cdot), ~ 
\hat{R}^{ji}_{nm}:= \br (u^i_n,u^j_m).
\end{eqnarray*}
The Hopf dual of the Hopf algebra $\A$ is denoted by $\A^\circ$. 
The L-functionals ${l^{\pm}}^i_j$ belong to $\A^\circ$. From (\ref{Kr0}) 
it follows that 
\begin{eqnarray*}
\Delta({l^{\pm}}^i_j) = \sum_{k=1}^n {l^{\pm}}^i_k \otimes {l^{\pm}}^k_j, 
~ i,j=1,{\cdots},n.
\end{eqnarray*}
These and the following 
relations will be often used in this paper:
\begin{eqnarray*}
(l^{+i}_{~j},u^k_l)=\hat{R}^{ik}_{lj},~~ (l^{-i}_{~j},u^k_l)
=(\hat{R}^{-1})^{ik}_{lj}=\bar{\br}(u^i_j,u^k_l),\\
(S(l^{+i}_{~j},u^k_l)=(\hat{R} ^{-1})^{ki}_{jl},~~
(S(l^{-i}_{~j},u^k_l)=\hat{R}^{ki}_{jl}. 
\end{eqnarray*}
Formula (\ref{Kr1}) implies that the matrix 
$\hat{R}$ and hence also $\hat{R}^{-1}$ intertwine the tensor 
product corepresentation $\uu \otimes \uu$.

Suppose that $\X$ is a subalgebra $\X$ of $\A$ such that 
$\Delta(\X) \subseteq \A \otimes \X$. Then $\X$ is a 
left $\A$-comdodule algebra or equivalently  
a {\it left quantum space} of $\A$ with left coaction  
$\varphi$ given by the restriction $\Delta \lceil\X$ 
of the comultiplication of $\A$. As in [KS],
such a subalgebra $\X$ will be called 
a {\it left quantum homogenous space} of 
the Hopf algebra $\A$. 

A {\it first order differential calculus} (abbreviated, a FODC) over
$\X$ is an $\X$-bimodule $\Gamma$ equipped with a linear mapping
$\kd:\X\rightarrow\Gamma$, called the differentiation, such
that:

\begin{description}
\item[(i)] $\kd$ satisfies the Leibniz rule $\kd(xy)=x{\cdot}\kd y+
\kd x{\cdot} y$ for any $x,y\in\X$,
\item[(ii)] $\Gamma$ is the linear span of elements $x{\cdot}\kd y{\cdot} z$
with $x,y,z\in\X$.
\end{description}

A FODC $\Gamma$ over $\X$ is called {\it left-covariant} if there exists a
linear mapping $\Phi:\Gamma\rightarrow\X\otimes\Gamma$  such that
$\Phi(x\kd y)=\Delta(x)(\id\otimes\kd)\Delta(y) $
for all $x,y \in\X$. For a left-covariant FODC $\Gamma$ of $\X$ 
the elements of the vector space 
$\Gammalinv =\{\eta\in\Gamma\,|\,\Phi(\eta)=1\otimes\eta\}$
are called {\it left-invariant} one-forms. A left-covariant FODC 
$\Gamma$ of $\X$ is called {\it inner} if there exists a
left-invariant one-form $\theta\in\Gammalinv$ such that
$$\kd x=\theta x-x\theta,~ x\in\X.$$

Let $\Gamma$ be a left-covariant FODC on the Hopf algebra $\A$ itself such 
that ${\rm dim}~ \Gamma := {\rm dim}~ \Gammalinv$ is finite-dimensional. 
We briefly recall a few facts from the general theory of these calculi 
(see [W2], [AS] or [KS], Section 14.1, for more details) that will be used 
in what follows. Such a FODC $\Gamma$ is characterized by a finite-dimensional 
subspace $\T$ of $\A^\circ$, called the {\it quantum tangent space} 
of $\Gamma$, and there is a canonical non-generate bilinear form 
$(.,.)$ on $\T \times \Gammalinv$.  If $\{X_i;i \in I \}$ and 
$\{\theta_i;i \in I \}$ are dual bases of $\T$ and $\Gammalinv$ with respect 
to this bilinear form, then the differentiation $\kd$ of the FODC 
$\Gamma$ can be expressed by 
\begin{eqnarray}\llabel{k2}
\kd a = 
\sum\nolimits_i a_{(1)}X_i(a_{(2)}) \theta_i,~ a \in \A.
\end{eqnarray}
The commutation 
relations between the elements of $\A$ and left-invariant one-forms 
of $\Gamma$ are given by 
\begin{eqnarray}\llabel{k3}
\theta_i a = \sum\nolimits_k a_{(1)}f^i_k(a_{(2)}) \theta_k , ~a \in \A,
\end{eqnarray} 
where $f^i_k$ are the functionals on $\A$ are determined by the 
equation 
\begin{eqnarray}\llabel{k4}
\Delta (X_k) - \varepsilon \otimes X_k = 
\sum\nolimits_i X_i \otimes f^i_k.
\end{eqnarray}
Let $\omega :\A \to \Gammalinv$ be the canonical projection 
defined by $\omega (a) 
= S(a_{(1)}) \kd a_{(2)}$ for $a \in \A $. Then one has
\begin{eqnarray}\llabel{bl}
(X,\omega (a)) = (X,a)~ \mbox{for} ~ X \in \T ~ \mbox{and}~ a \in \A.
\end{eqnarray}

If $\Gamma$ is a FODC of $\A$ with differentiation $\kd$, 
then $\tilde{\Gamma} := \X{\cdot} \kd \X {\cdot} \X$ is obviously a FODC of the 
subalgebra $\X$ with differentiation $\kd \lceil \X$. We call 
$\tilde{\Gamma}$ the {\it induced} FODC of the FODC $\Gamma$ of $\A$. Clearly,
if $\Gamma$ is left-covariant on the quantum group $\A$, then so is 
$\tilde{\Gamma}$ on the left quantum space $\X$. 

Our constructions of 
left-covariant FODC on $\A$ are based on the following lemma.
\mn

{\bf Lemma 1.} {\it A finite-dimensional vector space $\T$ of 
$\A^\circ$ is the quantum tangent space of a left-covariant 
FODC of $\A$ if and only if $X(1)=0$ and 
$\Delta(X)-\varepsilon\otimes X\in\X\otimes\A^\circ$ 
for all $X\in\X$.}
\mn

{\bf Proof.} [SS2], Lemma 1, or [KS], Proposition 14.5.\hfill $\Box{}$

\section*{3\ Quantum spaces generated by a row of ${\bf u}$}
Let $\X$ denote the unital subalgebra of $\A$ generated by the entries of the
last row of the matrix $\uu$, that is, by the elements $x_i:=u^i_n,i=1,{\dots},n$.
Clearly, $\X$ is a left quantum homogeneous space of $\A$ with left coaction
$\varphi=\Delta \lceil\X$ determined by
\begin{eqnarray}\llabel{Bus1}
\varphi(x_i)\equiv\Delta (u^i_n)={\sum\limits^n_{j=1}} u^i_j\otimes x_j,
~i=1,{\dots},n.
\end{eqnarray}
In this section we shall construct an $n$-dimensional 
left-covariant FODC $\Gamma$
on the Hopf algebra $\A$ which induces a FODC $\Gamma^\X$ on $\X$ such that
the differentials ${\kd}x_1,{\dots},{\kd}x_n$ form a free left $\X$-module basis
of $\Gamma^\X$.

First we define a FODC $\Gamma$ of $\A$. Let $\T^\X$ be the linear
span of functionals
$$
X_i:=\alpha^{-1} l^{-n}_{~~i} l^{-n}_{~~n},i=1,{\dots},n-1, ~{\rm and}~ X_n:=\alpha^{-1}
((l^{-n}_{~~n})^2-\varepsilon)
$$
on $\A$, where $\alpha$ is a non-zero complex number that will be specified
by formula (\ref{Bus5}) below. We assume that
\begin{equation}\llabel{Bus2}
l^{-m}_{~~n}=0~{\rm if}~m<n.
\end{equation}
Since $\Delta (l^{-n}_{~~n})=\sum_il^{-n}_{~~i}\otimes l^{-i}_{~~n}$, this assumption
implies in particular that $\Delta (l^{-n}_{~~n})=l^{-n}_{~~n}\otimes 
l^{-n}_{~~n}$, so that $l^{-n}_{~~n}$ is a character of the algebra $\A$ (that is,
$l^{-n}_{~~n} (ab)=l^{-n}_{~~n} (a) l^{-n}_{~~n}$ (b) for $a,b\in\A$ and $l^{-n}_{~~n}(1)=1)$.
Using the relation $\Delta (l^{-n}_{~~n})=l^{-n}_{~~n}\otimes l^{-n}_{~~n}$ 
we get
\begin{eqnarray*}
\Delta(X_i)-\varepsilon \otimes 
X_i&=&{\sum\limits_{j=1}^n}X_j\otimes l^{-j}_{~~i} l^{-n}_{~~n}, 
\, i=1,{\dots}, n-1,\\
\Delta (X_n)-\varepsilon \otimes X_n&=&X_n\otimes (l^{-n}_{~~n})^2.
\end{eqnarray*}
Because of (\ref{Bus2}), the latter equations can be written in the compact form
\begin{equation}\llabel{Bus3}
\Delta(X_i)-\varepsilon\otimes X_i={\sum\limits^n_{j=1}} X_i\otimes 
l^{-j}_{~~i} l^{-n}_{~~n},\, i=1,{\dots},n.
\end{equation}
Since obviously $X(1)=0$ and $\Delta (X)-\varepsilon\otimes X\in\T^\X
\otimes \A^0$  for all $X\in\T^\X$ by (\ref{Bus3}), 
it follows from Lemma 1 there exists a left-covariant
FODC $\Gamma$ on $\A$ such that $\T^\X$ is the quantum tangent space
of $\Gamma$.

Let us suppose in addition that
\begin{equation}\llabel{Bus4}
(l^{-n}_{~~i}, u^j_n)=0~{\rm if}~i\ne j~ ,i,j=1,{\dots},n,
\end{equation}
\begin{equation}\llabel{Bus5}
\alpha := (l^{-n}_{~~i}l^{-n}_{~~n}, u^i_n)=((l^{-n}_{~~n})^2, u^n_n)-1\ne 0 
~~{\rm for}~~ i=1,{\dots},n-1.
\end{equation}
 We abbreviate $c_-:=(l^{-n}_{~~n},u^n_n)$. Then we have $\alpha=c_-^{~~2}-1$. 

For $i=1,{\dots},n$, let $\theta_i$ denote the left-invariant
1-form $\omega(u^i_n)\equiv\sum_k S(u^i_k){\kd}u^k_n$ of $\Gamma$. The
assumptions (\ref{Bus4}) and (\ref{Bus5}) imply that 
$(X_j,u^i_n)=\delta_{ij}$ and so by formula (\ref{bl}) that
\begin{equation}\llabel{Bus6}
(X_j,\theta_i)=(X_j,\omega (u^k_n))=(X_j,u^i_n)=\delta_{ij}
\end{equation} 
for $i,j=1,{\dots},n$. In particular we conclude that the functionals
$X_1,{\dots},X_n$ are linearly independent, so that the FODC $\Gamma$
is $n$-dimensional. Further, (\ref{Bus6}) shows that $\{\theta_1,{\dots},\theta_n\}$
and $\{X_1,{\dots},X_n\}$ are dual bases of ${_{\rm inv}(\Gamma})$ and
$\T^X$, respectively. Therefore, comparing (\ref{k4}) and (\ref{Bus3}) 
and using (\ref{k3}), (\ref{Bus1}) and 
(\ref{Bus2}), we obtain for $r,j=1,{\dots},n$,
\begin{eqnarray}
\theta_rx_j&=& {\sum\limits_{k,s}} u^j_k (l^{-r}_{~~s} l^{-n}_{~~n}, u^k_n)\theta_s=
{\sum\limits_{k,m,s}} u^j_k (l^{-r}_{~~s}, u^k_m)(l^{-n}_{~~n},u^m_n)\theta_s\nonumber\\
\llabel{Bus7}
&=&{\sum\limits_{k,s}} c^{-1} (\hat{R}^{-1})^{rk}_{ns}~ u^j_k \theta_s.
\end{eqnarray}
These relations lead to the following commutation rules between the one-forms
$\theta_r$ and elements of the algebra $\X$:
\begin{equation}\llabel{Bus8}
\theta_r x={\sum\limits^n_{s=1}} x_{(1)} {\bar{\br}}(u^r_s, x_{(2)})(l^{-n}_{~~n},
x_{(3)})\, \theta_s,~x\in\X.
\end{equation}
Indeed, if $x$ is the generator $x_j$ of $\X$, then the third expression of
(\ref{Bus7}) can be rewritten as the right-hand side of (\ref{Bus8}). Using
the facts that $l^{-n}_{~~n}$ is a character and that ${\bar{\br}}_{21}$ is also
a universal $r$-form of $\A$ (see [KS], Proposition 10.2(iv)), one easily verifies that (\ref{Bus8}) holds for
a product $x^\prime x^{\prime\prime}$ provided that it holds for both
factors $x^\prime$ and $x^{\prime\prime}$. Thus, (\ref{Bus8}) is valid for
arbitrary elements $x$ of $\X$.\sn

Next we turn to the FODC $\Gamma^\X$ of $\X$.

\mn
{\bf Proposition 2.} (i) {\it The FODC $\Gamma$ of $\A$ induces a left-covariant
FODC $\Gamma^\X$ of $\X$ such that the set $\{{\kd} x_1,{\dots},{\kd} x_n\}$ is a
free left $\X$-module basis of $\Gamma^\X$. The $\X$-bimodule structure of 
$\Gamma^\X$ is determined by the commutation relations
\begin{equation}\llabel{Bus9}
{\kd} x_i{\cdot}x_j=(l^{-n}_{~~n}, u^n_n){\sum\limits^n_{k,m=1}} (\hat{R}^{-1})^{ij}_{km}
x_k{\cdot}{\kd} x_m, i,j=1,{\dots},n
\end{equation}
or equivalently by
\begin{equation}\llabel{Bus10}
{\kd} x_i{\cdot} x ={\sum\limits^n_{m=1}}\br (u^i_m, x_{(1)})x_{(2)} (l^{-n}_{~~n}, x_{(3)})
{\kd} x_m,\, x\in\X.
\end{equation}
(ii) For the differentiation ${\kd}$ of the FODC $\Gamma^\X$ of $\X$ we have}
\begin{equation}\llabel{Bus12}
{\kd} x = \alpha^{-1}(\theta_n x-x\theta_n),~~ x\in\X.
\end{equation}

\mn
{\bf Proof.} (i): First we prove formula (\ref{Bus9}). Since $(X_r, u^k_n)=
\delta_{kr}$, it follows from (\ref{k2}) and we have
\begin{equation}\llabel{Bus13}
{\kd} x_i\equiv{\kd} u^i_n ={\sum_{k,r}} u^i_k X_r (u^k_n)\theta_r=
\sum_r u^i_r\theta_r.
\end{equation}
Using (\ref{Bus7}), (\ref{Bus13}) and the fact that $\hat{R}^{-1}$ 
intertwines 
the tensor product corepresentation $\uu \otimes \uu$, we obtain
\begin{eqnarray*}
{\kd} x_i{\cdot}x_j&=&{\sum\nolimits_k} u^i_k\theta_k u^j_n={\sum\limits_{k,m,s}} c_- u^i_k u^j_m
(\hat{R}^{-1})^{km}_{ns}\theta_s\\
&=&{\sum\limits_{k,m,s}} c_-(\hat{R}^{-1})^{ij}_{km} u^k_n u^m_s\theta_s
={\sum\limits_{k,m}} c_-(\hat{R}^{-1})^{ij}_{km} x_k{\cdot}{\kd} x_m;
\end{eqnarray*}
which proves (\ref{Bus9}). Formula (\ref{Bus10}) can be derived from (\ref{Bus9})
similarly as (\ref{Bus8}) was from (\ref{Bus7}).\\
From (\ref{Bus9}) combined with the Leibniz rule it follows that $\Gamma^\X\equiv
\X{\cdot}{\kd}\X{\cdot}\X$ is equal to Lin$\{x\, {\kd} x_i; \, x\in\X,i=1,{\dots},n\}.$
Suppose that $\sum_ia_i{\kd} x_i=0$ for certain elements $a_i\in\X$. Then we
have $\sum_{i,k} a_iu^i_k\theta_k=0$. Since $\{\theta_1,{\dots},\theta_n\}$
is a free left $\A$-module basis of  $\Gamma^\X$, the latter yields $\sum_i
a_iu^i_k=0$ for $k=1,{\dots},n$ and hence $\sum_{i,k,j} a_iu^i_k S(u^k_j)=a_j=0$ 
for all $j=1,{\dots}, n$. Thus, $\{{\kd} x_1,{\dots},{\kd} x_n\}$ is a free left
$\X$-module basis of $\Gamma^\X$.\\
(ii): By (\ref{Bus4}) and (\ref{Bus5}) we have $c_-(\hat{R}^{-1})^{nk}_{ns} = 
(l^{-n}_{~~s}, u^k_n)(l^{-n}_{~~n},u^n_n) = 
\delta_{ks}(l^{-n}_{~~s}(l^{-n}_{~~n},u^k_n) = \delta_{ks} \alpha$ and 
$c_-(\hat{R}^{-1})^{nk}_{nn} = c_-(l^{-n}_{~~n},u^k_n) 
= \delta_{kn} c^{~2}_-$ 
for $s=1,\dots,n-1$ and $k=1,\dots,n.$ Inserting this into (\ref{Bus7}) and using (\ref{Bus13}) we 
obtain
\begin{eqnarray*}
\theta_n x_j&=& \sum\limits_{k=1}^{n-1} \alpha u^j_k \theta_k + 
c_-^{~2} u^j_n \theta_n\\
&=&\sum\limits_{k=1}^n \alpha u^j_k \theta_k + (c_-^{~2}-\alpha) u^j_n \theta_n 
=\alpha {\kd} x_j + x_j \theta_n,
\end{eqnarray*}
which proves (\ref{Bus12}) in the case $x=x_j$. Since both sides 
of (\ref{Bus12}),
considered as mappings of $\X$ to $\Gamma^\X$, satisfy the Leibniz rule,
(\ref{Bus12}) holds for all $x\in\X$.\hfill $\Box{}$

\mn
{\bf Remarks:}~ 1.) Since the left-invariant form $\theta_n\in\Gamma$ does not
belong to the $\X$-bimodule $\Gamma^\X$, formula (\ref{Bus12}) does
not mean that the FODC $\Gamma^\X$ is inner. It 
expresses rather the differentiation ${\kd}$ of $\Gamma^\X$ by means of an
extended bimodule in the sense of Woronowicz (see [W1]). But for the FODC
$\Gamma_1^\Z$ of the larger algebra $\Z$ considered in Section 4
the form $\theta_n$ is in $\Gamma_1^\Z$ and makes $\Gamma_1^\Z$ into
an inner FODC (see Proposition 4(iii) below).

2.) If $\A$ is one of the coordinate Hopf algebras 
$\cO(G_q), G_q=GL_q(N), SL_q(N), O_q(N),\break Sp_q(N)$, then the 
conditions (\ref{Bus2}) and (\ref{Camp1}) below can be assumed without loss 
of generality. This follows from the particular form of the universal R-matrix for 
the corresponding Drinfeld-Jimbo algebras 
(see, for instance, [KS], Theorem 8.17).

\bn
\section*{3\ Quantum spaces generated by a row of $\uu^{{\bf c}}$}
Lee $\Y$ be the subalgebra of $\A$ generated by the elements $y_i:=
(\uu^c)^i_n\equiv S(u^n_i),i=1,{\dots},n$, of the last row of the contragredient
corepresentation $\uu^c$. Then $\Y$ is a left quantum space of $\A$ with
left coaction $\varphi=\Delta\lceil\Y$ given on the generators
$y_i$ by 
$$
\varphi(y_i)\equiv\Delta (S(u^n_i))=
{\sum\limits^n_{j=1}} S(u^i_j)\otimes
y_j, \, i=1,{\dots},n.
$$
We shall proceed in a similar manner as in the preceding section. But the
considerations are technically slightly more complicated, because
we have to deal with square and inverse of the antipode of $\A$. 

Let $\beta$ be a non-zero complex number and let $\T^\Y$ be the linear span of
functionals 
$$
Y_i:=\beta^{-1} S(l^{+i}_{~~n}) l^{-n}_{~~n}, 
i=1,{\dots},n-1,~ {\rm and}~
Y_n:=\beta^{-1}((l^{-n}_{~~n})^2-\varepsilon).
$$
We assume that
\begin{equation}\llabel{Camp1}
l^{+n}_{~~m} =l^{-m}_{~~n}=0\mbox{ if }m<n\mbox{ and } S(l^{\pm n}_{~~n})=
l^{\mp n}_{~~n}.
\end{equation}
Similarly as in the preceding section, we then get 
\begin{equation}\llabel{Camp2}
\Delta(Y_i)-\varepsilon\otimes Y_i={\sum\limits^n_{j=1}} Y_j\otimes 
S(l^{+i}_{~~j})l^{-n}_{~~n},\, i=1,{\dots},n
\end{equation}
and $\T^\Y$ is the quantum tangent space of a left-covariant FODC
$\Gamma$ on $\A$.

Let us suppose in addition that there are numbers $\gamma_i\ne 0, i=1,{\dots},
n$, such that
\begin{equation}\llabel{Camp3}
S^2(u^i_j)=\gamma_iu^i_j\gamma_j^{-1},\, i,j=1,{\dots},n,
\end{equation}
and that
\begin{equation}\llabel{Camp4}
(l^{+i}_{~~n},u^n_j)=0\mbox ~{\rm for }~ i \ne j,\, i,j=1,{\dots},n,
\end{equation}
\begin{equation}\llabel{Camp5}
\beta := (l^{+n}_{~~n}l^{+i}_{~~n}, u^n_i)=
((l^{+n}_{~~n})^2, u^n_n)-1\ne 0 ~{\rm for}~ i=1,{\dots},n-1.
\end{equation}
We set $c:=(l^{+n}_{~~n},u^n_n)$ and $\eta_i 
:=\omega(S^{-1}(u^n_i))={\sum\nolimits_j} u^j_i{\kd} S
^{-1}(u^n_j)$ for 
$i=1,{\dots},n$. Since $S(l^{\pm n}_{~~n})=l^{\mp n}_{~~n}$ 
by (\ref{Camp1}), we have $c_-=(l^{-n}_{~~n},u^n_n)=c^{-1}$ and $\beta =c^2 -1 $.
 It is straightforward to check that (\ref{Camp1}), (\ref{Camp4})
and (\ref{Camp5}) imply that
\begin{equation}\llabel{Camp6}
(Y_j,\eta_i)=(Y_j, S^{-1} (u^n_i))=\delta_{ij}\mbox{ for }i,j=1,{\dots},n.
\end{equation}
Therefore, the FODC $\Gamma^\Y$ is $n$-dimensional. From (\ref{Camp2}),
(\ref{Camp1}) and (\ref{Camp3}) we get
\begin{eqnarray}
\eta_r y_j&=&{\sum\limits_{k,m,s}} S(u^k_j)(S(l^{+s}_{~~r}),S(u^m_k))
(l^{-n}_{~~n}, S(u^n_m))\eta_s\nonumber\\
\llabel{Camp7}
&=&{\sum\limits_{k,s}} c\gamma_n\gamma^{-1}_k {\hat{R}}^{sn}_{kr} 
S(u^k_j)\eta_s.
\end{eqnarray}
for $j,r=1,{\dots},n$. The first equality combined with the formulas 
$(S(l^{+s}_{~~r}), \cdot)=\br(S(\cdot),u^s_r)=\bar{\br}(\cdot,u^s_r)$ 
leads to the following form of the commutation relations
$$
\eta_ry={\sum\limits^n_{s=1}} y_{(1)}\bar{\br}(y_{(2)}, 
u^s_r)(l^{-n}_{~~n},y_{(3)})\eta_s,\, y\in\Y~.
$$
Let $\Gamma^\Y:=\Y{\cdot}{\kd}\Y{\cdot}\Y$ be the FODC on $\Y$ induced by the FODC
$\Gamma$ on $\A$.

\mn
{\bf Proposition 3.} (i) {\it $\Gamma^\Y$ is a left-covariant FODC on $\Y$ with the 
free left $\Y$-module basis $\{{\kd} {y_1},{\dots},{\kd} {y_n}\}$ and with $\Y$-bimodule
structure given by the relations
\begin{equation}\llabel{Camp8}
{\kd}{y_i}{\cdot} y_{j}=(l^{+n}_{~~n}, u^n_n){\sum\limits^n_{k,m=1}}\hat{R}^{mk}_{ji} 
y_k{\cdot} {\kd} y_m,~ i,j=1,
{\dots},n,
\end{equation}
or equivalently by
\begin{equation}\llabel{Camp9}
{\kd} y_i{\cdot} y={\sum\limits^n_{m=1}}\bar{\r} 
(y_{(1)}, u^m_i) y_{(2)}(l^{-n}_{~~n},
y_{(3)}) {\kd} y_m,~~
y\in\Y~.
\end{equation}
(ii) For any  $y\in\Y$ we have $ {\kd} y = \beta^{-1} (\eta_n y -y \eta_n).$}
\mn

{\bf Proof.} (i): It suffices to prove formula (\ref{Camp8}). First we note that (\ref{Camp6}) and (\ref{Camp3}) imply that
\begin{eqnarray}
{\kd}y_i&=&{\sum\limits_{k,r}} S(u^k_j)(Y_r,S(u^n_i))\eta_r={\sum\limits_{k,r}}
\gamma_n\gamma^{-1}_r S(u^k_i)(Y_r,S^{-1} (u^n_k))\eta_r\nonumber\\
\llabel{Camp12}
&=&\sum_r\gamma_u\gamma_r^{-1} S(u^r_i)\eta_r.
\end{eqnarray}
If $\br$ is a universal $r$-form of $\A$, then so is $\bar{\br}_{21}$ and
we have $\bar{\br}(a,S(b))=\br(a,b)$ and $\br (S(a),b)=\bar{\br} (a,b)$, 
where $\bar{\br}_{21} (a,b):=\bar{\br} (b,a)$ and $a,b\in\A$ 
(see, for instance,  
[KS]). Using these facts and the formulas (\ref{Kr1}) 
applied to $\bar{\br}_{21}$, (\ref{Camp7}), (\ref{Camp4}) and (\ref{Camp3}),
we compute
\begin{eqnarray*}
{\kd} y_i{\cdot} y_j&=& 
{\sum\limits_r} \gamma_n\gamma^{-1}_r S(u^r_i)\eta_r y_j\\
&=&{\sum\limits_{k,m,r,s}}\gamma_n\gamma^{-1}_r S(u^r_j) 
S(u^k_j)(l^{+l}_{~~n},
u^n_n)(l^{+s}_{~~r},S^2(u^n_k))u^m_s {\kd}S^{-1} (u^n_m)\\
&=&{\sum\limits_{k,m}} c\gamma_n\gamma^{-1}_r S(u^r_i)\left( 
{\sum\limits_{k,s}}S(u^k_j)u^m_s \bar{\br}_{21}(u^s_r, S(u^n_k))\right)
{\kd} S^{-1} (u^n_m)\\
&=&{\sum\limits_{k,m}} c\gamma_n\gamma^{-1}_i S^{-1}(u^r_i)
\left({\sum\limits_{k,s}}u^s_r S(u^n_k)\bar{\br}_{21}(u^m_s, S(u^k_j))\right) 
{\kd} S^{-1} (u^n_m)\\
&=&{\sum_{k,m}} c S(u^n_k)\bar{\br}(S(u^k_j),\gamma_m u^m_i\gamma^{-1}_i)
{\kd} S^{-1} (\gamma_n u^n_m\gamma^{-1}_m)\\
&=&{\sum\limits_{k,m}} cS(u^n_k)\bar{\br}(S(u^k_j),S^2(u^m_i)){\kd} 
S(u^n_m)\\
&=&{\sum\limits_{k,m}} c\hat{R}^{mk}_{ji} y_k{\cdot}{\kd} y_m.
\end{eqnarray*}
(ii): Since $c\hat{R}^{sn}_{kn} =(l^{+n}_{~~n},u^n_n)(l^{+s}_{~~n}, u^n_k) = 
\delta_{ks} (l^{+n}_{~~n}l^{+k}_{~~n},u^n_k) = \delta_{ks} \beta$ and 
$c\hat{R}^{nn}_{kn} = \delta_{kn} c^2$ by (\ref{Camp4}) and (\ref{Camp5}) 
for $s=1,\dots,n-1$ and $k=1,\dots,n$, it follows from (\ref{Camp7}) and 
(\ref{Camp12}) that
\begin{eqnarray*}
\eta_n y_j &=& \sum\limits^{n-1}_{k=1} \gamma_n \gamma_k^{-1} 
\beta S(u^k_j) \eta_k + c^2 S(u^n_j) \eta_n\\
&=& \sum\limits^n_{k=1}  \gamma_n \gamma_k^{-1} 
\beta S(u^k_j) \eta_k + (c^2- \beta) S(u^n_j) \eta_n
 = \beta {\kd} y_j + y_j \eta_n
\end{eqnarray*}
which implies the assertion.\hfill $\Box{}$
\bn

\section*{4\ Quantum spaces generated by a row of ${\bf u}$ and of 
${\bf u^c}$}
Let $\Z$ denote the subalgebra of $\A$ generated by the elements $x_i=
u^i_n$ and $y_i=S(u^n_i), i=1,{\dots},n$. That is, $\Z$ is the subalgebra of $\A$
generated by the algebras $\X$ and $\Y$. Our aim in this section is to 
construct four classes $\Gamma_j^\Z$, $j=1,2,3,4$, of left-covariant 
FODC  of $\Z$. 

First let us fix some notations and assumptions which will be kept 
in force throughout the whole section. Let $Z_n$ be a fixed 
group-like element of $\A^\circ$, that is, $Z_n(1)=1$
and $\Delta (Z_n)=Z_n\otimes Z_n$. Then $Z_n$ is invertible in $\A^\circ$
with inverse $Z_n^{-1}=S(Z_n)$. We retain the assumptions (\ref{Bus4}), 
(\ref{Camp1}), (\ref{Camp3}) and (\ref{Camp4}). In addition we suppose that 
\begin{equation}\llabel{Dur10}
(S(l^{+i}_{~~n}), u^j_n)=(l^{-n}_{~~i}, S^{-1}(u^n_j))=0\mbox{ for }(i,j)\ne
(n,n),i,j=1,{\dots},n,
\end{equation}
\begin{equation}\llabel{Dur11}
\gamma :=(l^{-n}_{~~i}~ l^{+n}_{~~n}, u^i_n)\ne 0\mbox{ and }\zeta :=(
l^{-n}_{~~n}~ l^{+i}_{~~n},u^n_i)\ne 0\mbox{ are independent}\\
\mbox{ of } i=1,{\dots},n-1,
\end{equation}
\begin{equation}\llabel{Dur12}
(l^{+n}_{~~n}, u^j_n)=(l^{-n}_{~~n},u^u_j)=(Z_n,u^j_n)=(Z^{-1}_{~~n},u^n_j)=0
\mbox{ if } i\ne n,
\end{equation}
\begin{eqnarray}\llabel{Dur13}
\delta := (Z_n,u^n_n)\ne 1.
\end{eqnarray}
Clearly, we then have $\delta^{-1}=(Z^{-1}_n, u^n_n)$. 

We now begin with the construction of the FODC $\Gamma_1^\Z$. Let $\T^\Z_1$ 
denote the
linear span of functionals
\begin{eqnarray*}
&&X_i:=\gamma^{-1}\delta^{-1} ~ l^{-n}_{~~i}~l^{+n}_{~~n}Z_n\mbox{ and } 
Y_i:=\zeta^{-1} \delta ~ S(l^{+i}_{~~n})l^{+n}_{~~n}Z_n,~ i=1,{\dots},n-1,\\
&&X_n :=(\delta-1)^{-1}(Z_n-\varepsilon)\mbox{ and }Y_n :=-\delta X_n=(\delta^{-1}
-1)^{-1} (Z_n-\varepsilon).
\end{eqnarray*}
For $i=1,{\dots},n-1$, we have
\begin{equation}\llabel{Dur14}
\Delta (X_i)-\varepsilon\otimes X_i=
\sum^{n-1}_{j=1} X_j\otimes l^{-j}_{~~i} 
l^{+n}_{~~n}Z_n+X_n\otimes (\delta-1)X_i,
\end{equation}
\begin{equation}\llabel{Dur15}
\Delta (Y_i)-\varepsilon\otimes Y_i=\sum^{n-1}_{j=1} Y_j\otimes 
S(l^{+i}_{~~j})~l^{-n}_{~~n} Z_n+Y_n\otimes (\delta^{-1}-1)Y_i,
\end{equation}
\begin{equation}\llabel{Dur16}
\Delta (X_n)-\varepsilon\otimes X_n=X_n\otimes Z_n,~~ \Delta(Y_n)-\varepsilon
\otimes Y_n=Y_n\otimes Z_n.
\end{equation}
Therefore, by Lemma 1, there exists a left-covariant
FODC $\Gamma_1$ on $\A$ with quantum tangent space $\T^\Z_1$. As in Sections 2
 and 3, we set 
$$\theta_j:=\omega(u^j_n)~ {\rm and}~ \eta_j:=\omega(S^{-1}(u^n_j)),~ 
j=1,{\dots},n,
$$
for the FODC $\Gamma^\Z$. From the assumptions (\ref{Bus4}), (\ref{Camp4}),
(\ref{Dur10}), (\ref{Dur12}) and the definition of the functionals $X_i,Y_i$ we immediately
derive
\begin{equation}\llabel{Dur17}
(X_r,u^s_n)=(Y_r, S^{-1}(u^n_s))=\delta_{rs}~ \mbox{ and} ~ 
(X_i,S^{-1}(u^n_j))=(Y_i, u^i_n)=0
\end{equation}
and so
$$
(X_r,\theta_s)=(Y_r,\eta_s)=\delta_{ks}\mbox{ and }(X_i,\eta_j)=(Y_i,\theta_j)=0
$$
for all $i,j,r,s=1,{\dots},n$ such that $(i,j)\ne (n,n)$. That is,
$\{\theta_1,{\dots},\theta_n,\eta_1,{\dots},\eta_{n-1}\}$ and $\{X_1,{\dots},
X_n, Y_1,{\dots},Y_{n-1}\}$ and likewise $\{\theta_1,{\dots},\theta_{n-1},
\eta_1,{\dots},\eta_n\}$ and $\{X_1,{\dots},X_{n-1},Y_1,{\dots},$\break $Y_n\}$
are dual bases of ${_{\rm inv}\Gamma}_1$ and $\T^\Z_1$, respectively. In particular,
we see that the FODC
 $\Gamma_1$ has the dimension ${\rm dim}~ \T^\Z_1=2n{-}1$. Moreover, the latter
facts imply that formula (\ref{Bus13}) and (\ref{Camp12}) hold for the differentiation 
${\kd}$ of the FODC $\Gamma_1$ as well. Further, from the formulas (\ref{k3}) and 
(\ref{Dur14})--(\ref{Dur16}) we obtain the following commutation 
relations between the basis elements of ${_{\rm inv}\Gamma}_1$ and elements $a\in\A$:
\begin{eqnarray}\llabel{Dur18}
\theta_ra&=&\sum^{n-1}_{s=1} a_{(1)} (l^{-r}_{~~s}l^{+n}_{~~n} 
Z_n,a_{(2)})\theta_s,~~
\eta_sa=\sum^{n-1}_{s=1}a_{(1)}
(S(l^{+s}_{~~r})l^{+n}_{~~n} Z_n,a_{(2)})\eta_s,\\
\llabel{Dur19}
\theta_na&=&a_{(1)}(Z_n,a_{(2)})\theta_n+
(\delta-1)\sum^{u-1}_{s=1} a_{(1)}((X_s,a_{(2)})
\theta_s+(Y_s,a_{(2)})\eta_s),\\
\llabel{Dur20}
\eta_n a&=& a_{(1)}(Z_n,a_{(2)})\eta_n+(\delta^{-1}-1)
\sum^{n-1}_{s=1} a_{(1)} ((X_s,a_{(2)})
\theta_s+(Y_s,a_{(2)})\eta_s)
\end{eqnarray}
for $r=1,{\dots},n-1$.

Let $\Gamma_1^\Z$ denote the FODC of $\Z$ which induced by the FODC 
$\Gamma_1$ of $\A$.

\bn
{\bf Proposition 4.} (i) {\it For the $\Z$-bimodule $\Gamma_1^\Z$ we have 
the commutation relations
\begin{eqnarray*}
{\kd} x_i{\cdot} x_j&=&c\delta\sum^n_{k,m=1}
 (\hat{R}^{-1})^{ij}_{km} 
x_k{\cdot} {\kd} x_m + (\delta - \gamma\delta -1)(x_i{\cdot} {\kd} x_j-
x_ix_j\theta_n),\\
{\kd} y_i{\cdot} y_j&=&c^{-1}\delta^{-1} 
\sum^n_{k,m=1}\hat{R}^{mk}_{ji}y_k{\cdot}
{\kd} y_m+ (\delta^{-1}- \zeta \delta^{-1}-1)(y_i{\cdot} {\kd} y_j- 
y_iy_j\eta_n),\\
{\kd} x_i{\cdot} y_j&=&c^{-1}\delta^{-1} 
\sum^n_{k,m=1}\hat{R}^{ki}_{mj}y_k{\cdot}
{\kd} x_m+(\delta-1)(x_i{\cdot}{\kd} y_j-x_iy_j\eta_n),\\
{\kd} y_i{\cdot} x_j&=& c\delta\sum^n_{k,m=1}(\grave{R}^-)^{ij}_{km} x_k{\cdot}{\kd} 
y_m+(\delta^{-1}-1)(y_i{\cdot}{\kd} x_j-y_ix_j\theta_n),
\end{eqnarray*}
where $(\grave{R}^-)^{ij}_{km}:=\bar{\br}(u^j_k,S^2(u^m_i)),
~ i,j,k,m =1,{\dots},n$. \\ 
(ii) The set $\{\kd x_1,{\dots},\kd x_n,\kd y_1,{\dots},\kd y_n\}$ generates
$\Gamma_1^\Z$ as a left $\Z$-module. For arbitrary 
elements $a_1,{\dots},a_n,b_1,{\dots},b_n\in\Z$, the relation
\begin{equation}\llabel{Dur4}
{\sum\limits^n_{i=1}} (a_i\kd x_i +b_i\kd y_i)=0
\end{equation}
is equivalent to the following set of equations:
\begin{equation}\llabel{Dur5}
a_j=\left( {\sum\limits^n_{i=1}} a_ix_i \right) y_j,~~ b_j=\left( 
{\sum\limits^n_{i=1}}b_ix_i\right)
x_j\gamma_j\gamma^{-1}_n\mbox{ for }j=1,{\dots},n,
\end{equation}
\begin{equation}\llabel{Dur6}
{\sum\limits^n_{i=1}} a_ix_i=(l^{+n}_{~~n},u^n_n)~{\sum\limits^n_{i=1}} b_i
y_i.
\end{equation}
(iii) $\Gamma_1^\Z$ is an inner FODC of $\Z$ with respect to the left-invariant
one-form $\theta_n=-\delta \eta_n$, that is, we have}
\begin{equation}\llabel{Dur21}
{\kd} z=(\delta-1)^{-1} (\theta_nz-z\theta_n) ~~ for~~ z\in\Z.
\end{equation}

\mn
{\bf Proof.} (i): We carry out the proofs of the second and the fourth
relations and work with the dual bases $\{ \theta_1, {\dots}, \theta_{n-1},
\eta_1, {\dots}, \eta_n \}$ and $\{ X_1, {\dots}, 
X_{n-1}, Y_1, {\dots}, Y_n\}$.
The two other relations follow by a similar slightly simpler reasoning.
Using the formulas (\ref{Dur16}), (\ref{Dur18}), (\ref{Dur20}) and the above 
assumptions we compute
\begin{eqnarray*}
{\kd} y_i{\cdot} y_j&=&\sum_r \gamma_n\gamma_r^{-1}S(u^r_i)\eta_r S(u^n_j)\\
&=&\sum^{n-1}_{r=1} \sum^n_{k,s=1} \gamma_n\gamma^{-1}_r S(u^r_i) S(u^k_j)
S(l^{+s}_{~~r})l^{+n}_{~~n}Z_n,S(u^n_k))\eta_s\\
&&+\sum^{n-1}_{s=1}\sum^n_{k=1} S(u^n_i) S
(u^k_j)(\delta^{-1}-1)((X_s,S(u^n_k))
\theta_s+(Y_s,S(u^n_k))\eta_s)\\
&&+\sum^n_{k=1}S(u^n_i)S(u^k_j)(Z_n,S(u^n_k))\eta_n\\
&=&\sum^n_{k,r,s=1}\gamma_n\gamma_r^{-1} S(u^r_i)S(u^k_j)\delta^{-1}c^{-1}
(l^{+s}_{~~r},S^2(u^n_k))\eta_s\\
&&+\sum^{n-1}_{s=1}\sum^n_{k=1} S(u^u_i S(u^k_j)
(\delta^{-1}-1-\zeta\delta^{-1})(Y_s,S(u^n_k))\eta_s.
\end{eqnarray*}
The first sum is treated as in the proof of Proposition 3. In this 
manner it becomes equal to
$c^{-1}\delta^{-1}\sum_{k,m}\hat{R}^{mk}_{ji}y_k{\kd} y_m.$ 
Put $\tilde{\zeta}:=\zeta \delta^{-1}+1-\delta^{-1}$. 
Since $(Y_s,S(u^n_k))=\gamma_n\gamma^{-1}_k(Y_s,
S^{-1}(u^n_k))=\gamma_n\gamma^{-1}_k \delta_{ks}$ by (\ref{Dur15}) and 
$\gamma_n\gamma^{-1}_k \eta_k=\gamma_n\gamma^{-1}_k 
\omega(S^{-1}(u^n_k))=\omega(S(u^n_k))=\sum_r S^2(u^r_k){\kd} 
S(u^n_r)$, the 
second expression yields 
\begin{eqnarray*}
\quad
&&\sum^{n-1}_{k=1} -\tilde{\zeta} S(u^n_i)S
(u^k_j)\gamma_n\gamma^{-1}_k\eta_k\\
&&\qquad\qquad=\tilde{\zeta} S(u^n_i)S
(u^n_j)\eta_n-\sum^n_{k,r=1} \tilde{\zeta}S(u^n_i)
S(u^k_j)S^2(u^r_k){\kd} S(u^n_r)\qquad\qquad\\
&&\qquad\qquad=\tilde{\zeta} y_iy_j\eta_n-\tilde{\zeta} y_i{\cdot}{\kd} y_i.\qquad\qquad
\end{eqnarray*}
Putting both terms together we obtain the second relation. In order to prove
the fourth relation we proceed in a similar manner. Using the facts that
$(S(l^{-s}_{~~n})l^{+n}_{~~n} Z_n,u^k_n)=\zeta\delta^{-1} (Y_s,u^k_n)=0$
and $(X_s,u^k_n)=\delta_{ks}$ for $s=1,{\dots},n-1$, we obtain
\begin{eqnarray*}
{\kd} y_i{\cdot} x_j=&&\sum^{n-1}_{r=1} \sum^n_{k,s=1} \gamma_n\gamma^{-1}_r S(u^r_i)u^j_k
(S(l^{+s}_{~~r}) l^{+n}_{~~n} Z_n,u^k_n)\eta_s\\
&&+\sum^{n-1}_{s=1} \sum^n_{k=1} S(u^n_i)u^j_k(\delta^{-1}-1)((X_s,u^k_n)
\theta_s+(Y_s,u^k_n)\eta_s)\\
&&+\sum^n_{k=1} S(u^n_i)u^j_k (Z_n,u^k_n)\eta_n\\
=&&\sum^n_{k,m,r,s=1} c\delta\gamma_n\gamma^{-1}_r S(u^r_i)u^j_ku^m_s
(\hat{R}^{-1})^{ks}_{rn} {\kd} S^{-1}(u^n_m)\\
&&+\sum^{n-1}_{k=1} (\delta^{-1}-1)S(u^n_i)u^j_k\theta_k\\
=&&\sum^n_{k,m,r,s=1} c\delta\gamma_n\gamma^{-1}_r S(u^r_i)u^k_ru^s_n\bar{\br}
(u^j_s,u^m_k){\kd} S^{-1}(u^n_m)\\
&&-(\delta^{-1}-1)S(u^n_i)u^j_n\theta_n+\sum^n_{k=1}(\delta^{-1}-1)S
(u^n_i)u^j_k\theta_k\\
=&&\sum^n_{k,m,r,s=1} c\delta S(u^r_i)S^2
(u^k_r)u^s_n \bar{\bf r} (u^j_s,S^2(u^m_k))
{\kd} S^{-1}(S^2(u^n_m))\\
&&-(\delta^{-1}-1)y_ix_j\theta_n+(\delta^{-1}-1)y_i{\cdot}{\kd} x_j\\
=&&\sum^n_{m,s=1} c\delta (\grave{R}^-)^{ij}_{sm} x_s {\cdot} {\kd} y_m 
+(\delta^{-1}-1) (y_i{\cdot}{\kd} x_j-y_i x_j\theta_n).
\end{eqnarray*}
(ii): Since  $\theta_n=\sum_i y_i{\kd} x_i$ and $\eta_n=\sum_i\gamma_n\gamma^{-1}_i
x_i{\kd} y_i,$ the four relations in (i) imply that the set $\{{\kd} x_1,{\dots},
{\kd} x_n, {\kd} y_1,{\dots},{\kd} y_n\}$ generates $\Gamma_1^\Z$ 
as a left $\Z$-module. It remains to verify that (\ref{Dur4}) is 
equivalent to (\ref{Dur5}) and (\ref{Dur6}). 
Since $X_n=-\delta Y_n$, the element $S^{-1}(u^n_n)+\delta u^n_n$ is 
annihilated by the whole quantum tangent space $\T_1^\Z$ and hence 
$0=\omega(S^{-1}(u^n_n)+\delta u^n_n)=\eta_n+\delta \theta_n$. 
Inserting the relations $\eta_n=-\delta \theta_n$, (\ref{Bus13}) 
and (\ref{Camp12}) into (\ref{Dur4}) we see that (\ref{Dur4}) reads as 
$$
{\sum\limits^n_{i=1}}\left( 
{\sum\limits^{n-1}_{r=1}}(a_iu^i_r \theta_r+b_i\gamma_n
\gamma_r^{-1} S(u^r_i)\eta_r)+(a_iu^i_n-c^2 b_i S(u^n_i))\theta_n
\right) =0.$$
Since the set $\{ \theta_1,{\dots},\theta_n,\eta_1,{\dots},\eta_{n-1} \}$
is a free left $\A$-module basis of $\Gamma_1^{\Z}$, the latter is equivalent
to the relations
\begin{eqnarray}\llabel{Dur8}
&& {\sum\nolimits_i} a_iu^i_r={\sum\nolimits_i} b_i S(u^r_i)=0 \mbox{ for }
r=1,{\dots}, n-1,\\
\llabel{Dur9}
&& {\sum\nolimits_i} (a_iu^i_n - \delta b_i S(u^n_i))=0.
\end{eqnarray}
Multiplying $\sum_i a_iu^i_r$ by $S(u^r_k)$ and $\sum_i b_i S(u^r_i)$ by 
$S^2(u^k_r)=\gamma_k\gamma^{-1}_r u^k_r$ and summing over $r$, (\ref{Dur8})
implies (\ref{Dur5}). Formula (\ref{Dur9}) is nothing but (\ref{Dur6}). Using
the relations $\sum_i y_ix_i=\sum_ix_iy_i\gamma_i\gamma_n^{-1}=1$,
equation (\ref{Dur5}) in turn implies (\ref{Dur8}).\\
(iii): It suffices to prove (\ref{Dur21}) for the generators
$z=x_i,y_i$. Because $(Y_s, u^j_n)=0$ and $(X_s, u^j_n)=\delta_{js}$
for $s=1,{\dots},n-1$, it follows from (\ref{Dur19}) and (\ref{Bus13}) that
\begin{eqnarray*}
\theta_nx_i&=&u^i_n (Z_n,u^n_n)\theta_n+(\delta-1)\sum^{n-1}_{s=1} u^i_s
\theta_s\\
&=&u^i_n(\delta-(\delta-1))\theta_n+(\delta-1)\sum^n_{s=1} u^i_s\theta_s\\
&=&x_i\theta_n+(\delta-1){\kd} x_i,
\end{eqnarray*}
which gives (\ref{Dur21}) in the case $z=x_i$. Similarly, 
using the formulas (\ref{Dur20}) and
(\ref{Camp12}) we get ${\kd} y_i=(\delta^{-1}-1)^{-1}(\eta_n Y_i-Y_i\eta_n),$ 
so that ${\kd} y_i=(\delta-1)^{-1}(\theta_ny_i-y_i\theta_n)$.\hfill $\Box{}$
\bn

Next we turn to the FODC $\Gamma_4^\Z$ on $\Z$. We take the linear span 
$\T^\Z_4$ of functionals 
\begin{eqnarray*}
X_i:=\gamma^{-1} l^{-n}_{~~i}l^{+n}_{~~n} ~{\rm and} ~ 
Y_i:=\zeta^{-1} S(l^{+1}_{~~n})l^{+n}_{~~n},~ i=1,{\dots},n-1,\\
X_n=(\delta-1) ^{-1}(Z_n-\varepsilon) ~{\rm and} ~ 
Y_n=(\delta^{-1}-1)^{-1} (Z_n-\varepsilon).
\end{eqnarray*}
For $i=1,{\dots},n-1$, we then have 
\begin{equation}\llabel{Dur45}
\Delta (X_i)-\varepsilon\otimes X_i={\sum\limits^{n\cdot 1}_{j=1}} X_j
\otimes l^{-j}_{~~i} l^{+n}_{~~n},~~ \Delta(X_n)-\varepsilon\otimes X_n=X_n
\otimes Z_n,
\end{equation}
\begin{equation}\llabel{Dur46}
\Delta (Y_i)-\varepsilon\otimes Y_i={\sum\limits^{n-1}_{j=1}} 
Y_j\otimes S(l^{-i}_{~~j}) l^{+n}_{~~n},~~\Delta (Y_n)-\varepsilon\otimes Y_n=Y_n\otimes 
Z_n.
\end{equation}
These formulas and the relations $(X_i,u^n_j)=(Y_i,S^{-1}(u^n_j))=
\delta_{ij}$ and $(X_i, S(u^n_k))=(Y_i, u^k_n)=0$ for $i,j=1,{\dots},n$
and $k=1,{\dots},n-1$ imply that $\T^\Z_4$ is the quantum tangent space
of a $(2n{-}1)$-dimensional left-covariant FODC $\Gamma_4$ of $\A$. The
commutation relations of this FODC between the one-forms $\theta_r,\eta_s$ and 
elements of $\A$ are
\begin{eqnarray*}
\theta_ra&=&\sum^{n-1}_{s=1}
 a_{(1)}(l^{-r}_{~~s} l^{+n}_{~~n},a_{(2})\theta_s,~~
\eta_ra=\sum^{n-1}_{s=1} a_{(1)}((S(l^{+s}_{~~r})l^{+n}_{~~n},
a_{(2)})\eta_s,\\
\theta_n a&=&a_{(1)}(Z_n,a_{(2)})\theta_n,~~ \eta_n a=
a_{(1)}(Z_n,a_{(2)})\eta_n
\end{eqnarray*}
for $r,s=1,{\dots},n-1$. Let $\Gamma_4^\Z$ denote the FODC of $\Z$ 
which is induced by the FODC $\Gamma_4$ of $\A$. By similar computations as
carried out above one prove the following cross commutation relations of
the $\Z$-bimodule $\Gamma_4^\Z$:
\begin{eqnarray*}
{\kd} x_i{\cdot} x_j&=& c\sum^n_{k,m=1} (\hat{R}^{-1})^{ij}_{km} x_k
{\cdot}{\kd} x_m - \gamma x_i{\cdot}{\kd} x_j+ \gamma x_ix_j\theta_n,\\
{\kd} y_i{\cdot} y_j&=&c^{-1}\sum^n_{k,m=1}\hat{R}^{mk}_{ji} y_k{\cdot} 
{\kd} y_m -\zeta y_i{\cdot} {\kd} y_j+ \zeta y_iy_j\eta_n,\\
{\kd} x_i{\cdot} y_j&=&c^{-1}\sum^n_{k,m=1}\hat{R}^{ki}_{mj}~ y_k{\cdot}
{\kd} x_m,\\
{\kd} y_i{\cdot} x_j&=& c\sum^n_{k,m=1} (\grave{R}^-)^{ij}_{km} 
x_k{\cdot}{\kd} y_m.
\end{eqnarray*}
These are precisely the relations which are obtained 
by setting formally $\delta = 1$ in the commutation relations 
for the FODC $\Gamma_1^\Z$ (see Proposition 4(i)). That is, the FODC 
$\Gamma_4^\Z$ can be viewed as the limit of the FODC $\Gamma_1^\Z$ as 
$\delta \to 1$. Note that the FODC $\Gamma_1^\Z$ has no direct meaning in the case 
$\delta = 1$.

By "mixing" the elements of the quantum tangent spaces of the FODC $\Gamma_1^\Z$ and $\Gamma_4^\Z$ 
one obtains two other FODC $\Gamma_1^\Z$ and $\Gamma_4^\Z$ on $\Z$. We briefly 
describe the quantum Lie algebras of the corresponding FODC 
$\Gamma_2$ and $\Gamma_3$ of $\A$ and the 
commutation rules of these calculi. Let $\T^\Z_2$ be the linear span of 
functionals 
\begin{eqnarray*}
&&X_i:=\gamma^{-1}\delta^{-1}~ l^{-n}_{~~i}~l^{+n}_{~~n}Z_n\mbox{ and } 
Y_i:=\zeta^{-1} ~ S(l^{+i}_{~~n})l^{+n}_{~~n},~ i=1,{\dots},n-1,\\
&&X_n :=(\delta-1)^{-1}(Z_n-\varepsilon)\mbox{ and }Y_n :=-\delta X_n=(\delta^{-1}
-1)^{-1} (Z_n-\varepsilon).
\end{eqnarray*}
and $\T^\Z_3$ the span of functionals
\begin{eqnarray*}
&&X_i:=\gamma^{-1} ~ l^{-n}_{~~i}~l^{+n}_{~~n}\mbox{ and } 
Y_i:=\zeta^{-1} \delta ~ S(l^{+i}_{~~n})l^{+n}_{~~n}Z_n,~ i=1,{\dots},n-1,\\
&&X_n :=(\delta-1)^{-1}(Z_n-\varepsilon)\mbox{ and }Y_n :=-\delta X_n=(\delta^{-1}
-1)^{-1} (Z_n-\varepsilon).
\end{eqnarray*}
From the formulas (\ref{Dur14}),(\ref{Dur15}),(\ref{Dur45}) and 
(\ref{Dur46}) we see that $\T^\Z_2$ and $\T^\Z_3$ are quantum tangent spaces of $(2n{-}1)$-dimensional 
left-covariant FODC $\Gamma_2$ and $\Gamma_3$ of $\A$, respectively. 
From these formulas we also read off the following commutation relations between the 
left-invariant one-forms $\theta_i,\eta_k$ and elements $a\in\A$:
\begin{eqnarray*}
\Gamma_2:~~ 
\theta_ra&=&\sum^{n-1}_{s=1} a_{(1)} (l^{-r}_{~~s}l^{+n}_{~~n} 
Z_n,a_{(2)})\theta_s,~~
\eta_sa=\sum^{n-1}_{s=1}a_{(1)}
(S(l^{+s}_{~~r})l^{+n}_{~~n},a_{(2)})\eta_s,\\
\theta_na&=&a_{(1)}(Z_n,a_{(2)})\theta_n+
(\delta-1)\sum^{u-1}_{s=1} a_{(1)}(X_s,a_{(2)})
\theta_s,\\
\eta_n a&=& a_{(1)}(Z_n,a_{(2)})\eta_n+(\delta^{-1}-1)
\sum^{n-1}_{s=1} a_{(1)} ((X_s,a_{(2)})\theta_s,\\
\Gamma_3:~~ \theta_ra&=&\sum^{n-1}_{s=1} a_{(1)} (l^{-r}_{~~s}
l^{+n}_{~~n},a_{(2)})\theta_s,~~
\eta_sa=\sum^{n-1}_{s=1}a_{(1)}
(S(l^{+s}_{~~r})l^{+n}_{~~n} Z_n,a_{(2)})\eta_s,\\
\theta_na&=&a_{(1)}(Z_n,a_{(2)})\theta_n+
(\delta-1)\sum^{n-1}_{s=1} a_{(1)}(Y_s,a_{(2)})\eta_s,\\ 
\eta_n a&=& a_{(1)}(Z_n,a_{(2)})\eta_n+(\delta^{-1}-1)
\sum^{n-1}_{s=1} a_{(1)}(Y_s,a_{(2)})\eta_s,
\end{eqnarray*}
where $r=1,{\dots},n-1$. As earlier, the FODC
on $\Z$ induced by the FODC $\Gamma_j$ on $\A$ is denoted
by $\Gamma_j^\Z$, $j=2,3$. From the preceding set of formulas one gets 
the following commutation rules for the $\Z$-bimodule $\Gamma_j^\Z$: 
\begin{eqnarray*}
\Gamma_2^\Z:~~ {\kd} x_i{\cdot} x_j&=&c\delta\sum^n_{k,m=1}
 (\hat{R}^{-1})^{ij}_{km} 
x_k{\cdot} {\kd} x_m + (\delta - \gamma\delta -1)(x_i{\cdot} {\kd} x_j-
x_ix_j\theta_n),\\
{\kd} y_i{\cdot} y_j&=&c^{-1} \sum^n_{k,m=1}\hat{R}^{mk}_{ji}y_k{\cdot}
{\kd} y_m - \zeta y_i{\cdot} {\kd} y_j + \zeta y_iy_j\eta_n,\\
{\kd} x_i{\cdot} y_j&=&c^{-1}\delta^{-1} 
\sum^n_{k,m=1}\hat{R}^{ki}_{mj}y_k{\cdot}{\kd} x_m,\\ 
{\kd} y_i{\cdot} x_j&=& c\sum^n_{k,m=1}(\grave{R}^-)^{ij}_{km} x_k{\cdot}{\kd} 
y_m+(\delta^{-1}-1)y_i{\cdot}{\kd} x_j\\
\Gamma_3^\Z: ~~{\kd} x_i{\cdot} x_j&=&c\sum^n_{k,m=1}
 (\hat{R}^{-1})^{ij}_{km} 
x_k{\cdot} {\kd} x_m - \gamma x_i{\cdot} {\kd} x_j + \gamma 
x_ix_j\theta_n,\\
{\kd} y_i{\cdot} y_j&=&c^{-1}\delta^{-1} 
\sum^n_{k,m=1}\hat{R}^{mk}_{ji}y_k{\cdot}
{\kd} y_m+ (\delta^{-1}- \zeta \delta^{-1}-1)(y_i{\cdot} {\kd} y_j- 
y_iy_j\eta_n),\\
{\kd} x_i{\cdot} y_j&=&c^{-1} 
\sum^n_{k,m=1}\hat{R}^{ki}_{mj}y_k{\cdot}
{\kd} x_m+(\delta-1)x_i{\cdot}{\kd} y_j,\\
{\kd} y_i{\cdot} x_j&=& c\delta\sum^n_{k,m=1}(\grave{R}^-)^{ij}_{km} 
x_k{\cdot}{\kd} y_m.
\end{eqnarray*}
Recall that by Proposition 4(iii) the FODC $\Gamma_1^\Z$ of $\Z$ is inner. It 
turns out that none of the three other FODC 
$\Gamma_2^\Z$, $\Gamma_3^\Z$, $\Gamma_4^\Z$ 
is inner. Indeed, from the
above commutation rules one easily derives that
\begin{equation}\llabel{Dur22}
\Gamma_2^\Z:~~ \theta_n y_i=\delta^{-1} y_i\theta_n, ~~
\Gamma_3^\Z:~~ \theta_n x_i=x_i\theta_n ,~~
\Gamma_4^\Z:~~ \theta_n y_i= y_i \theta_n
\end{equation}
for all $i=1,{\dots},n$ . Further, for all four FODC $\Gamma_j^\Z$ 
we have $\theta_n=-\delta\eta_n$ and this is up to complex multiples 
the only left-invariant one-form
of $\Gamma_j^\Z$. Therefore, we conclude at once from (\ref{Dur22}) 
that none of the FODC $\Gamma_j^\Z$, $j=2,3,4,$ of $\Z$ is inner. 

All four left-covariant FODC $\Gamma_j^\Z$ of $\Z$ depend on 
the group-like element
$Z_n\in\A^\circ$. It can be freely choosen such that it satisfies 
the conditions 
(\ref{Dur12}) and (\ref{Dur13}). This dependence is
reflected by the appearance of the parameter $\delta =(Z_n,u^n_n)$ 
in the above formulas. For the FODC $\Gamma_1^\Z$ a distinguished choice 
of $Z_n$ is $Z_n=(l^{-n}_{~~n})^2$. In this case $\T^\Z_1$
is just the sum of the quantum tangent spaces $\T^\X$ and $\T^\Y$ considered
in Sections 2 and 3 and the FODC~ $\Gamma_1^\Z$ might be thought as gluing 
together the FODC $\Gamma^\X$ and $\Gamma^\Y$. Further, if we assume 
in addition the conditions (\ref{Bus5}) and (\ref{Camp5}), then we have
$\alpha=\gamma c^{-2} =\delta-1,\beta=\gamma c^2=\delta^{-1}-1$, and $\delta
=c^{-2}$,  so that $\gamma\delta+1-\delta=\zeta\delta^{-1}+1-\delta^{-1}=0$. 
Thus, in this case the first two relations for the FODC $\Gamma_1^\Z$ 
in Proposition 4(i) become even linear.
\bn

\section*{5.\ Application to the quantum homogeneous space ${\bf GL_q(N)
/ GL_q(N{-}1)}$}
In this section let $\A$ denote the Hopf algebra $\cO(GL_q(N)), 
\uu=(u^i_j)_{i,j=1,{\dots},N}$ the fundamental corepresentation 
of $\A$ and ${\hat R}$
the corresponding $R$-matrix given by (see [FRT])
\begin{equation}\llabel{elb1}
R^{ji}_{kl}\equiv
\hat{R}^{ij}_{kl}:=q^{\delta_{ij}} \delta_{il}\delta_{jk}+
(q-q^{-1})\theta (j-i) \delta_{ik} \delta_{j},~~ i,j,k,l=1,{\dots},N.
\end{equation}
The Hopf algebra $\A$ is coquasitriangular with universal $r$-form $\br$
determined by
\begin{equation}\llabel{elb2}
\br (u^i_j,u^k_l)={\hat R}^{ki}_{jl},~ i,j,k,l=1,{\dots},N.
\end{equation}
Further, We suppose that $Z_n$ is a monomial in the main diagonal 
L-functionals $l^{\pm i}_{~~j}$. 

Using (\ref{elb1}) and (\ref{elb2}) one easily verifies that the above 
assumptions (\ref{Bus2}), (\ref{Bus4}), (\ref{Bus5}), 
(\ref{Camp1}), (\ref{Camp3}), (\ref{Camp4}), 
(\ref{Camp5}), (\ref{Dur10}), (\ref{Dur11}), (\ref{Dur12}) and 
(\ref{Dur13}) are then fulfilled 
with $n=N$, $\alpha=-\zeta =q^{-2}-1,\beta=- \gamma =q^2-1,c=q$ and 
$\gamma_i=q^{2i}$. Therefore, 
{\it all results obtained in the preceding three sections are 
valid in this case}. Here we shall add only a few remarks concerning 
these results rather than 
restating them in the present situation. The quantum homogeneous 
space $\X$ is then, of course, isomorphic to 
the quantum vector space $\cO({\bbbc}^N_q)$ ([KS], Proposition 9.11) and the FODC $\Gamma^\X$ is 
one of the two well-known covariant calculi on $\cO({\bbbc}^N_q)$ 
discovered in [PW] and [WZ]. However, the approach given in Section 2 might be 
still of interest. The FODC $\Gamma_j^\Z, j=1,2,3,4,$ developed in the 
preceding section are left-covariant FODC on the subalgebra $\Z$ of $\A$ generated by the element $x_i\equiv u^i_N$
and $y_i\equiv S(u^N_i), i=1,{\dots},N$. All four FODC 
have the property that $\Gamma_j^\Z$ as a left $\Z$-module is generated by the
differentials ${\kd} x_1,{\dots},{\kd} x_N,{\kd} y_1,{\dots},{\kd} y_N$. 
The FODC $\Gamma_1^\Z$ described by Proposition 4 is inner. In the special case 
$Z_n=(l^{-n}_{~~n})^2$ it coincides with the distinguished calculus 
considered in [We] (more precisely with its left-convariant counter-part). 

The importance of the left quantum
space $\Z$ stems from the fact that it is (isomorphic to) the quantum homogeneous 
space $GL_q(N)/ GL_q(N{-}1)$. Indeed, there is a unique surjective 
Hopf algebra homomorphism $\pi:GL_q(N)\rightarrow GL_q(N{-}1)$ such that
\begin{eqnarray*}
&&\pi(u^i_j)=w^i_j,~ i,j=1,{\dots},N{-}1,\\
&&\pi (u^i_N)=\pi(u^N_i)=0,~ i=1,{\dots}, N{-}1,~ \pi (u^N_N)=1,
\end{eqnarray*}
where $w^i_j, i,j=1,{\dots}, N{-}1$, denote the matrix entries of the fundamental
matrix for the quantum group $GL_q(N{-}1)$. Then the set
$$
\cO(GL_q(N)/ GL_q(N{-}1)):=\{a\in\cO(GL_q(N)):({\rm id}\otimes \pi)
\circ\Delta (a)=a\otimes 1\}
$$
of all right $GL_q(N{-}1)$-invariant elements of $\cO(GL_q(N))$ is a subalgebra
and a left quantum space for $\cO(GL_q(N))$ with respect to the coaction $\Delta\lceil
\cO(GL_q(N))$. Obviously, the elements $x_i$ and $y_i$ are in $\cO(GL_q(N)
/ GL_q(N{-}1))$, so that $\Z\subseteq \cO(GL_q(N)/ GL_q(N{-}1))$.
If $q$ is not a root of unity, then we have the equality $\Z=\cO(GL_q(N)
/ GL_q(N{-}1))$. (For  the corresponding right quantum space $GL_q(N{-}1)
\backslash GL_q(N)$
this is proved in [NYM], Proposition 4.4, or [KS], Section 14.6. The proof 
for the left quantum space 
$GL_q(N)/ GL_q(N{-}1)$ is completely analogous.) 

 Suppose now that $q$ is a real number and $q\ne 0,\pm 1$. Then it is 
well-known that the Hopf algebra $\cO(GL_q(N))$ is a Hopf $\ast$-algebra,
denoted by $\cO(U_q(N))$, with involution determined by  $(u^i_j)^\ast=
S(u^j_i),i,j{=}1,{\dots},N$. Further, the algebra $\cO(GL_q(N)/ 
GL_q(N{-}1))$ is a $\ast$-subalgebra such that $x^\ast_i\equiv
(u^i_N)^\ast =y_i\equiv S(u^N_i)$ and a left $\ast$-quantum space for $\cO(U_q(N))$. It
is denoted by $\cO(U_q(N)/ U_q(N{-}1))$ and called the coordinate
$\ast$-algebra of the {\it quantum sphere} associated with the quantum group
$U_q(N)$. In this case the two left-covariant FODC $\Gamma_1$ and 
$\Gamma_4$ of $\cO(U_q(N))$ 
and hence their induced FODC $\Gamma_1^\Z$ and $\Gamma_4^\Z$ 
on $\Z=\cO(U_q(N)/ U_q(N{-}1))$ are $\ast$-calculi.
We prove these assertions for $\Gamma_1$ and $\Gamma_1^\Z$. First note that 
$(l^{\pm i}_{~~j})^\ast
=S(l^{\pm j}_{~~i})$ (see [KS], formula (10.47)) for the corresponding involution
of the Hopf dual $\cO(GL_q(N))^\circ$. Hence we obtain $X^\ast_N=X_N$ and
$X^\ast_i=(l^{-N}_{~~i}
l^{-N}_{~~N}Z_N)^\ast=Z_N l^{-N}_{~~N} S(l^{+i}_{~N})$
for $i=1,{\dots},N{-}1$. Since $Z_N$ is a monomial in the L-functionals 
$l^{\pm i}_{~~i}$, $Z_N l^{-N}_{~~N} S(l^{+i}_{~N})$ is a  complex multiple 
of $S(l^{+i}_{~N})l^{-N}_{~~N}Z_N =Y_i$. Therefore, we have $X^\ast\in\T_1^\Z$
for all $X\in\T_1^\Z$, so that
$\Gamma_1^\Z$ is a $\ast$-calculus of $\cO(U_q
(N))$ by Proposition 14.6 in [KS]. Since $\Z$ is a $\ast$-subalgebra of
$\cO(U_q(N))$, the induced FODC $\Gamma_1^\Z$ is also a $\ast$-calculus. 
Thus, {\it the FODC
$\Gamma_1^\Z$ and $\Gamma_4^\Z$ are $\ast$-calculi on the coordinate $\ast$-algebra $\Z=
\cO(U_q(N)/ U_q(N{-}1))$ of the quantum sphere}. Note that because
these FODC are $\ast$-calculi it suffices to prove only one of the 
commutation relations for ${\kd}x_i{\cdot}x_j$ and ${\kd}y_i{\cdot}y_j$ 
and one of the relations for ${\kd}x_i{\cdot}y_j$ and 
${\kd}y_i{\cdot}x_j$. The two others follow then by
applying the involution and inverting the corresponding $R$-matrix. The FODC 
$\Gamma_2^\Z$ and $\Gamma_3^\Z$ are not $\ast$-calculi on 
$\Z$, but one has $(\T_2^\Z)^\ast = \T_3^\Z$.
 
Let us return to the general case where $q$ is a complex number such that 
$q\ne 0,\pm 1$. From its very construction it is clear that the left-covariant 
$(2n{-}1)$-dimensional FODC $\Gamma_1$ of the 
Hopf algebra $\cO(GL_q(N))$ is a useful tool for the study of the 
induced FODC $\Gamma_1^\Z$ on the subalgebra
$\Z$. However, $\Gamma_1$ is not suitable as a FODC of the Hopf algebra
$\cO(GL_q(N))$ itself, because the generators $X_i,Y_i$ of the quantum tangent
space $\T^\Z$ are only supported on the last row and column of the 
fundamental 
matrix $\uu=(u^i_j)$. To remedy this defect, one can construct an $N^2$-dimensional 
left-covariant FODC $\Gamma$ on $\A=\cO(GL_q(N))$ that induces the FODC 
$\Gamma_1^\Z$ on $\Z$ as well. We restrict ourselves to the 
distinguished calculus $\Gamma_1^\Z$ with $Z_n=(l^{-n}_{~~n})^2$. 
Let $\T$ be the linear span of linear functionals
\begin{eqnarray}\llabel{elb3}
X_{ij}&=&(q^{-2}-1)^{-1}l^{-j}_{~~i} l^{-j}_{~~j},~ i<j,\\
\llabel{elb4}
Y_{ji}&=&(q^2-1)^{-1} S(l^{+i}_{~~j})l^{-j}_{~~j},~ i<j,\\
\llabel{elb5}
X_{ii}&=&(q^{-2}-1)^{-1}((l^{-i}_{~~i})^2-\varepsilon),~ 
Y_{ii}=-q^{-2}X_{ii}
\end{eqnarray}
on $\A$. For $i\le j, i,j=1,{\dots},N$, we then have
\begin{eqnarray}\llabel{elb6}
&&\Delta (X_{ij})-\varepsilon\otimes X_{ij}=\sum_{k\le j} X_{kj}\otimes 
l^{-k}_{~~i}l^{-j}_{~~j},\\
\llabel{elb7}
&&\Delta (Y_{ji})-\varepsilon\otimes Y_{ji}=\sum_{k\le j} Y_{jk}\otimes 
S(l^{+i}_{~~k}) l^{-j}_{~~j}.
\end{eqnarray}
Thus, by Lemma 1, there is a left-covariant FODC $\Gamma$
of $\A$ which has the quantum tangent space $\T$. From the explicit form 
(\ref{elb1}) of the matrix $\hat{R}$ and its inverse $\hat{R}^{-1}=(q-q^{-1})
\hat{R}+I$ we compute 
\begin{equation}\llabel{elb8}
(X_{ij},u^r_s)=\delta_{ir}\delta_{js}\mbox{ and }  (Y_{ji},  
S^{-1} (u^s_r))=\delta_{ir} \delta_{js}\mbox{ for }i\le j,
\end{equation}
Setting $\theta_{ij}:=\omega(u^i_j)={\sum\nolimits_k} S(u^i_k)du^k_j$ and 
$\eta_{ji}:=\omega (S^{-1}(u^j_i))={\sum\nolimits_k} u^k_i d S^{-1} (^j_k)
\mbox{ for } i\le j$, the formulas (\ref{bl}) and (\ref{elb8}) imply that
\begin{equation}\llabel{elb9}
(X_{ij}, \theta_{rs})=(Y_{ji},\eta_{sr})=\delta_{ir}\delta_{js}\mbox{ for }
i\le j, i,j,r,s=1,{\dots},N.
\end{equation}
In particular, the functionals $X_{ij}, Y_{rs}, i\le j, s<r$, are linearly 
independent, so that the FODC $\Gamma$ has dimension $N^2$. Further, it follows 
from (\ref{bl}) and (\ref{elb9}) that the sets $\{\theta_{ij},\eta_{rs}; i\le j, s<r\}$ and $\{X_{ij},
Y_{rs}; i\le j, s<r\}$ and also the sets $\{\theta_{ij},\eta_{rs}; i<j, s\le r\}$ 
and $\{ X_{ij},Y_{rs}; i<j,s\le r\}$ are dual bases of ${_{\rm inv}\Gamma}$ and
 $\T$, respectively. It is not difficult to verify that the two calculi $\Gamma$
and $\Gamma_1$ with $Z_n=(l^{-n}_{~~n})^2$ of $\A$ induce the same 
FODC $\Gamma_1^\Z$ on the quantum
space $\Z$. 

For $j=1,{\dots},N$, let $\T_j$ denote the linear span of 
functionals $X_{ij}$ and $Y_{ji}, i\le j$. From (\ref{elb6}) and 
(\ref{elb7}) we conclude that there is a $(2j{-}1)$-dimensional 
left-covariant FODC
$\Gamma^j$ on $\cO(GL_q(N))$ which has the quantum tangent space $\T_j$.
The FODC $\Gamma^N$ is nothing but the FODC 
$\Gamma_1$ developed in Section 4 (as always throughout this 
discussion, with $Z_n=(l^{-n}_{~~n})^2$). Since the linear quantum tangent space 
$\T$ is the direct sum of vector spaces $\T_1,{\dots},\T_N$, the FODC 
$\Gamma$ is the direct sum of FODC $\Gamma^1,{\dots},\Gamma^N$. These and 
other properties indicate that the FODC 
$\Gamma$i is a promising tool for the study of the interplay between 
the quantum group $GL_q(N)$ and the quantum homogeneous spaces $GL_q(j)
/ GL_q(j{-}1), j=2,{\dots},N$. The FODC $\Gamma$ is only left-covariant, 
but not bicovariant. However, because of its particular and simple structure 
the FODC $\Gamma$ might be even more important and useful for appliations 
and computations than the bicovariant 
calculi of the Hopf algebra $\cO(GL_q(N))$. We shall return to this matter 
in Section 7.

At the end of this section, let us briefly turn to the quantum group 
$SL_q(N)$. The Hopf algebra $\cO(SL_q(N))$ is also quasitriangular with 
universal $r$-form $\br$ such that 
\begin{equation}\llabel{elb10}
\br(u^i_j,u^k_l)=z {\hat R}^{ki}_{jl}, i,j,k,l=1,{\dots}, N,
\end{equation}
where ${\hat R}$ is given by (\ref{elb1}) and $z$ is a complex $N$-th root of $q^{-1}$.
Then the linear span of functionals $X_{ij}, Y_{ij}, i<j$, and $X_{rr}, r=2,
{\dots},N$, defined by (\ref{elb3})--(\ref{elb5}) is also the quantum tangent space of a $(N^2-1)
$-dimensional FODC on $\cO(SL_q(N))$. It should be emphasized that because of 
the appearance of the number $z$ in (\ref{elb10}) the equalities 
(\ref{elb8}) are no longer valid for $\cO(SL_q(N))$. Some $(N^2{-}1)$-dimensional 
left-covariant FODC on 
$\cO(SL_q(N))$ with reasonable properties have been 
constructed in [SS3]. This FODC is different from those in [SS3], 
but it is based on a similar idea.
\bn

\section*{6.\ A left-covariant FODC on ${\bf GL_q(N)/ GL_q(N{-}1)}$ induced from a
bicovariant FODC on ${\bf GL_q(N)}$}
In this  section we retain the notation of the preceding section. Let 
$\Gamma_{\rm bi}$ be the bicovariant FODC on $\A=\cO(GL_q(N))$ constructed by the 
bicovariant bimodule $(u^c\otimes u,L^+\otimes L^{-,c})$. (Details can be 
found, for instance, in [KS], Sections 14.5 and 14.6). Here we only 
need the two 
facts (see [KS], 14.6.3 and Example 14.8) that the set  $\{\omega_{ij}:=
\omega (u^i_j)=\sum_k S(u^i_k){\kd} u^k_j,i,j=1,{\dots},N\}$ is a basis of the 
vector space ${_{\rm inv}(\Gamma}_{\rm bi})$ of left-invariant one-forms of $\Gamma_{\rm bi}$ 
and that the commutation rules between the forms $\omega_{ij}$ and an element 
$a \in\A$ are given by
\begin{equation}\llabel{for1}
\omega_{ij} a=\sum_{r,s}~ a_{(1)} l^{+i}_{~~r} (a_{(2)})S(l^{-s}_{~~j})
(a_{(3)})\omega_{rs}.
\end{equation}

\bn
{\bf Proposition 5.} {\it The FODC $\Gamma_{\rm bi}$ induces a left-covariant FODC
$\Gamma^\Z$ on the quantum space $\Z$ such that
\begin{eqnarray*}
{\kd} x_i {\cdot} x_j= q\sum_{k,m} \hat{R}^{ij}_{km} x_k{\cdot}{\kd} x_m,\\
{\kd} y_i{\cdot} y_j=q^{-1} \sum_{k,m} (\hat{R}^{-1} )^{ji}_{mk} y_k{\cdot}
{\kd} y_m,\\
{\kd} x_i{\cdot} y_j=q^{-1} \sum_{k,m} ( \hat{R}^{-1})^{ki}_{mj}
y_k{\cdot} {\kd} x_m,\\
{\kd} y_i{\cdot} x_j=q \sum_{k,m} \grave{R}_{km}^{ij} x_k {\cdot}{\kd} y_m,
\end{eqnarray*}
where $\grave{R}^{ij}_{km} := {\bf r} (S^2(u^k_j), u^i_m),~ i,j,k,m=1,
{\dots},N$. Further, we have $\omega_{NN} x_i=q^2 x_i \omega_{NN}$ and 
$\omega_{NN} y_i=q^{-2} y_i\omega_{NN}$ for $i=1,{\dots},N$.}

\mn
{\bf Proof.}  We verify, for instance, the third commutation relation. 
From the explicit form (\ref{elb1}) of the matrix ${\hat R}$ 
it follows that $({\hat R}^{-1})^{sN}_{lN}= q^{-1} \delta_{sN}\delta_{lN}$ 
for $s,l=1,{\dots},N$. Using essentially this 
fact and formula (\ref{for1}) we compute
\begin{eqnarray*}
{\kd} x_i{\cdot} y_j&=& \sum_k u^i_k \omega_{kN} S(u^N_j)\\
&=&\sum_{k,m,l,r,s} u^i_k S(u^m_j)(l^{+k}_{~~r}, S(u^l_m))(S(l^{-s}_{~~N}),
S(u^N_l))\omega_{rs}\\
&=&\sum_{r,s,l} \left( \sum_{k,m} u^i_k S(u^m_j)\left(\hat{R}^{-1}
\right)^{lk}_{rm}\right) q^{2N-2l} ( \hat{R}^{-1})^{sN}_{lN} 
\omega_{rs}\\
&=&\sum_{r,s,l} \left(\sum_{k,m} S(u^l_k) u^m_r 
(\hat{R}^{-1})^{ki}_{mj}\right)q^{2N-2l} q^{-1} \delta_{sN}
\delta_{lN} \omega_{rs}\\
&=&\sum_{k,m,r} q^{-1} (\hat{R}^{-1})^{ki}_{mj} S(u^N_k) u^m_r 
\omega_{rN}\\
&=&\sum_{k,m} q^{-1} (\hat{R}^{-1})^{ki}_{mj} y_k{\cdot} {\kd} x_m.
\end{eqnarray*}
The other relations follow by similar reasonings as
above or as used ealier. We shall not carry out the details. \hfill $\Box{}$
\bn

The FODC $\Gamma^\Z$ is another good candidate of a reasonable 
differential calculus on the quantum homogeneous space $\Z$. It is 
a $\ast$-calculus if $q$ is real and the involution of $\Z$ is given by 
$x^\ast_i=y_i,i=1,{\dots},N$, because the FODC $\Gamma_{\rm bi}$ on $\cO (GL_q(N)$ 
is known to be a $\ast$-calculus with respect to the involution $(u^i_j)^\ast=
S(u^j_i), i,j=1,{\dots},N.$ But there is a striking 
difference between the two distinguished calculi 
$\Gamma^\Z$ and $\Gamma_1^\Z$: $\Gamma_1^\Z$ is inner,
while $\Gamma^\Z$ is not. In order to verify the latter, it suffices to note 
that $\omega_{NN}x_i - x_i\omega_{NN} =(q^2-1) x_i\omega_{NN}$ is 
obviously not a multiple of ${\kd} x_i$.
\bigskip\noindent

\section*{7.\ A recipe for the construction of left-covariant FODC}
The first order differential calculi on quantum homogeneous spaces
developed above are induced from left-covariant calculi on 
the quantum group. All  these left-covariant calculi on the corresponding Hopf 
algebra are built by the same simple recipe that will be elaborated more
explicitely in this section. As always, $\A$ is a coquasitriangular Hopf
algebra and $l^{\pm i}_{ ~j}$ are the $L$-functionals on $\A$ with respect
to a fixed corepresentation $\uu=(u^i_j)_{ij=1,{\dots},n}$ of $\A$. Throughout
this section we retain assumption (\ref{Camp1})).

Let $i,j\in\{1,{\dots},n\}$ be two indices such that $i\le j$ and let $Z$ be
a group-like element of $\A^\circ$. Define
\begin{eqnarray*}
X^+_r&=&l^{+i}_{~~r} l^{-i}_{~~i} Z,~r=i+1,{\dots}j,
\mbox{ and } X^+_i =Z{-}\varepsilon,\\
X^-_r&=&l^{-j}_{~~r} l^{+j}_{~~j} Z,~r=i,{\dots},j{-}1,
\mbox{ and } X^-_j =Z{-}\varepsilon,\\
Y^+_r&=&S(l^{+r}_{~~j})l^{+j}_{~~j}Z,~r=i,{\dots},j{-}1,
\mbox{ and } Y^+_j = Z{-}\varepsilon,\\
Y^-_r&=&S(l^{-r}_{~~i})l^{-i}_{~~i}Z,~r=i{+}1,{\dots},j,
\mbox{ and } Y^-_i = Z{-}\varepsilon,\\
\T^\pm_{ij}(Z)&=& \mbox{Lin}~\{X^\pm_r; i\le r\le j\},~ \T^\pm_{ji} (Z)
=\mbox{Lin}~
\{Y^\pm_r;i\le r\le j\}.
\end{eqnarray*}
Using (\ref{Camp1}) one easily verifies 
that each vector space $\T=\T^\pm_{ij}(Z),
\T^\pm_{ij}(Z)$ has the properties that $X(1)=0$ and $\Delta (X)-\varepsilon
\otimes X\in\T\otimes \A^\circ$ for all $X\in\T$. Hence, by Lemma 1 each space $\T^\pm_{ij}
(Z),\T^\pm_{ji}(Z)$ is the quantum tangent space of a
left-covariant FODC $\Gamma^\pm_{ij}, \Gamma^\pm_{ji}$ on $\A$. 
Let us call the first order calculi of the form $\Gamma^\pm_{ij},\Gamma^\pm_{ji}$
{\it elementary} FODC. All left-covariant FODC on $\A$ occuring in this
paper are direct sums of elementary FODC (with possible different group-like
elements $Z$!). By forming sums of elementary FODC one gets a large
supply of left-covariant FODC which have a very simple structure and are
easy to handle. FODC of this form have been introduced in [SS3]. Note that the
commutation rules of the elements of the elements of the quantum tangent spaces
obtained in this manner are not necessarily quadratically closed and that 
the dimensions of the spaces of higher forms may be different from the corresponding
classical dimensions (see [SS2] for such examples). 

For the group-like
elements $Z$ one may take, for instance, a monomial in the main diagonal
$L$-functionals $l^{\pm i}_{~~i},i=1,{\dots},n$. Interesting choices
of $Z$ are, of course, $Z=l^{\pm i}_{~~i}$ for $\T^{\pm}_{ij} (Z)$ and 
$Z=l^{\mp i}_{~~i}$ for $\T^\pm_{ji}(Z)$ or $Z=\varepsilon$ for all four
FODC. Let us illustrate this by simple examples and set
\begin{eqnarray*}
\T^+&=& {\sum\nolimits_i} \T^+_{in} (l^{+i}_{~~i})=\mbox{Lin}~ \{l^{+i}_{~~j}; i\le 
j, i,j=1,{\dots},n\},\\
\T^-&=& {\sum\nolimits_j} \T^-_{ij} (l^{-j}_{~~j})=\mbox{Lin}~ \{l^{-j}_{~~i}; i\le 
j, i,j=1,{\dots},n\},\\
\T_+&=& {\sum\nolimits_j} \T^+_{ij} (l^{-j}_{~~j})=\mbox{Lin}~ \{S(l^{+i}_{~~j}); i\le 
j, i,j=1,{\dots},n\},\\
\T_-&=& {\sum\nolimits_i} \T^-_{ni} (l^{+i}_{~~i})=\mbox{Lin}~ \{S(l^{-j}_{~~i}); i\le 
j, i,j=1,{\dots},n\},
\end{eqnarray*}
Then, $\T^+,\T^-,\T_+,\T_-,\T^+{+}\T_-$ and $\T^-{+}\T_+$ are
quantum tangent spaces of left-covariant FODC on $\A$.

Now we want to be more specific and suppose that $\A$ is one of the
Hopf algebras $\cO(G_q),G_q=GL_q(N),SL_q(N), O_q(N), Sp_q(N)$, and $\uu$ is the
fundamental corepresentation.

\mn
{\bf Case 1.} $\A=\cO(GL_q(N))$\\
Then the vector spaces $\T^+{+}\T_-$ and $\T^-{+}\T_+$ defined above are the quantum
tangent spaces of two $N^2$-dimensional left-covariant FODC on $\cO(GL_q(N))$.
It is easily seen that the commutation relations of the elements of both quantum
tangent spaces are quadratically closed. Further, it can be shown that the
dimensions of the spaces of $k$-forms for the associated universal higher
order differential calculi (see [KS], 14.3, for this notion) are 
$ {N^2 \choose k} $ as
in the classical case.

\mn
{\bf Case 2.} $\A=\cO(SL_q(N))$\\
In this case, $\T^+{+}\T_-$ and $\T^-{+}\T_+$ are also $N^2$-dimensional FODC on
$\cO(SL_q(N))$, but we are interested in FODC that have the classical group
dimension $N^2{-}1$. It is rather easy to construct such FODC: Let $\T_{\rm od}$
be the sum of $\T^+_{in}(\varepsilon),\T^-_{1i}(\varepsilon),i=1,{\dots},n,$
and let $\T_{\rm md}$ be the vector space spanned by $N{-}1$ of the $N$ functionals
$l^{+i}_{~~i}-\varepsilon$. Then, $\T=\T_{\rm od}+\T_{\rm md}$ is the quantum
tangent space of an $(N^2{-}1)$-dimensional FODC on $\cO(SL_q(N))$. This
first order calculus strongly resembles the ordinary differential calculus
on the Lie group $SL(N)$ in many aspects. But it has the disadvantage that the
commutation rules between elements of the quantum tangent space (for instance, 
$l^{+i}_{~N} l^{-i}_{~~i}$ and $l^{-N}_{~~j} l^{+N}_{~~N}$) do not close 
quadratically. $(N^2{-}1)$-dimensional FODC on $O(SL_q(N)$) that do not
have this defect have been constructed in [SS3]. However, using the
same idea as in [SS3], the quantum tangent space $\T$ can be modified by
multiplying the secondary diagonal elements such that commutation relations close
quadratically.

In order to be more precise, let $f_i$ and $g_i, i=1,{\dots},N$, be monomials
in the main diagonal $L$-functionals $l^{\pm j}_{~~j}$. Let $\T_{\rm od}$ be the
linear span of $X_{ij}:=l^{+i}_{~~j} l^{-i}_{~~i} f_i$ and $X_{ji}:=l^{-j}_{~~i}
l^{+j}_{~~j} g_j, i<j$, and let $\T_{\rm md}$ be an 
$(N{-}1)$-dimensional vector space generated by functionals of the form 
$f-\varepsilon$, where $f$ is a monomial in $l^{\pm i}_{~~i}, j=1,{\dots},N$. 
Suppose that $f_i,g_i\in\bbbc\varepsilon \oplus\T_{\rm md}$ for $i=1,{\dots},N$.
Then one easily verifies that $\T:=\T_{\rm od}+\T_{\rm md}$ is the quantum tangent
space of an $(N^2{-}1)$-dimensional FODC on $\cO(SL_q(N))$. Further, the
commutation relations for elements of $\T$ are quadratically closed if 
and only if $f^{-1}_i g_i(l^{+i}_{~~i})^2$ is independent of 
$i=1,{\dots},N$. (This assertion and the explicit form of commutation rules
can be derived from the relations $L^\pm_1 L^\pm_2 R=RL^\pm_2L^\pm_1$ and 
$L^-_1 L_2^+ R=RL^+_2 L^-_1$ using (\ref{elb1}). We omit the details.) These conditions
can be fulfilled as follows: Fix an index $k\in\{1,{\dots},N\}$ and set
$g_i=(l^{-i}_{~~i})^2(l^{+k}_{~~k})^2$ and $f_i=\varepsilon$ for $i=1,{\dots},N$.
Another possible choice is $f_i=(l^{+i}_{~~i})^2 (l^{-k}_{~~k})^2$ and 
$g_i=\varepsilon$ for $i=1,{\dots},N$.  These two special cases are in fact
the two FODC $\Gamma_1$ and $\Gamma_2$ constructed in [SS3].

In order to come into contact with the considerations in Sections 4 and 5,
we carry out the same consideration based in the generators $X^-_r,
Y^+_r$ rather than $X^+_r, X^-_r$. We suppose that the elements $f_i,g_i$
and the vector space $\T_{\rm md}$ satisfy the assumptions stated in the first
half of the preceding paragraph. Now let $\T_{\rm od}$ be the vector space
generated by the functionals $X_{ji}=l^{-i}_{~~j} l^{+i}_{~~i}f_i$ and $X_{ij}=
S(l^{+i}_{~~j})l^{+j}_{~~j}g_j, i<j$. Then $\T:=\T_{\rm od}+\T_{\rm md}$ is again
the quantum tangent space of an $(N^2-1)$-dimensional FODC on $\cO(SL_q(N))$.
The commutation relations for $\T$ close quadratically if and only if $f_ig_i(l^{+i}_{~~i})^2$
does not depend on $i=1,{\dots},N$.

\mn
{\bf Case 3:} $\A=\cO(O_q(N))$ and $\A=\cO(Sp_q(N))$\\ 
In this case the fundamental matrix $\uu$ fulfills 
the metric condition
\begin{eqnarray}
\uu C \uu^t C^{-1} = C \uu C^{-1} u = I
\end{eqnarray}
and the R-matrix is given by 
\begin{eqnarray}
\hat{R}^{ji}_{mn}=q^{\delta_{ij}-
\delta_{ij^\prime}} \delta_{im}
\delta_{jn}+(q-q^{-1})\theta (i-m)(\delta_{jm}\delta_{in}-
\epsilon C^j_iC^m_n),
\end{eqnarray}
where $i^\prime := n+1-i$, $\epsilon =1$ for $O_q(N)$, 
$\epsilon =-1$ for $Sp_q(N)$ and $C=(C^i_j)$ is the corresponding 
matrix of the metric (see [FRT] or [KS] for details). We shall essentially 
use the fact that $C^i_j = 0$ if 
$i \ne j^ \prime$.

Before we turn to the construction of the FODC, let us look for a 
moment at the "ordinary" first order calculus on the Lie 
groups $O(N)$ and $Sp(N)$. Then the matrix
\begin{eqnarray}
\theta = S(\uu){\kd}\uu = (\theta_{ij} \equiv 
\sum\nolimits_k S(u^i_k) {\kd} u^k_j))_{i,j=1,\dots,N}
\end{eqnarray}
satisfies the relation
\begin{eqnarray}\llabel{r1}
\theta = -C^{-1} \theta C,~\mbox{ that is,}~ \theta_{ji} = (C^{-1})^i_{i'} \theta_{j'}^{i'}  
C_j^{j'}~ \mbox{for}~ i,j=1,\dots, N.
\end{eqnarray}
We briefly sketch the proof of this well-known fact. Indeed, differentiating
the condition $\uu C^{-1} \uu = C^{-1}$, we obtain 
\begin{eqnarray}\label{r2}
{\kd} \uu^t C^{-1} \uu + \uu^t C^{-1} {\kd} \uu = 0.
\end{eqnarray}
From  $C \uu^t C^{-1} \uu = I$ we get $C^{-1} S(\uu) = \uu^t C^{-1}$ and so 
$\uu^t C^{-1} {\kd} \uu = C^{-1} \theta.$ For the metric $C$ of the Lie groups $O(N)$ 
and $Sp(N)$ we have $(C^{-1})^t = \epsilon C^{-1}$. Hence the relation 
$C^{-1}S(\uu) = \uu^t C^{-1}$ implies that $C^{-1}\uu = S(\uu)^t C^{-1}$. 
Because functions and forms commute (!) for the classical 
differential calculus, we can write 
${\kd} \uu^t C^{-1} \uu = ( S(\uu) {\kd}  \uu )^t C^{-1} = 
\theta^t C^{-1}.$ 
Inserting these expressions into (\ref{r2})  we obtain (\ref{r1}).

For the construction of the left-covariant FODC we shall restrict ourselves 
to the quantum group $O_q(N)$. In the case of $Sp_q(N)$ 
one has to omit the elements $X_{ii}$ supporting the secondary diagonal entries 
$u^i_{i^\prime}$ in order to be in accordance with the ordinary 
calculus on the classical group $Sp(N)$. The remaining parts 
are verbatim the same.

Let us abbreviate $I:=\{ (i,j): i' \leq j, i,j = 1,\dots,N \}$. Then 
the elements $u^i_j$ with $(i,j) \in I$ are precisely those 
entries of the matrix $\uu$ that are below or on the secondary diagonal. Now we define
\begin{eqnarray*}
X_{ji}:= l^{-j}_{~~i} l^{+j}_{~~j}Z_j ~~ {\rm and} ~~ 
X_{ij} := l^{+j'}_{~~i'} l^{+j'}_{~~j'}Z_{j'}~~ {\rm for }~~ j' \leq i < j , 
i,j=1,\dots, n,
\end{eqnarray*}
where $Z_j$ and $Z_{j'}$ are group-like elements of the Hopf dual $\cO (O_q(N))^\circ$. 
These functionals $X_{ji}, X_{ij}$ separate the elements $u^i_j$ 
such that $(i,j) \in I$ and $i \neq j$. 
In order to separate also the main diagonal entries $u^i_{i'}, i' \leq i,$ we 
choose group-like elements $Y_i, i'  \leq i,$ of $\cO (O_q(N))^\circ$ such that 
\begin{eqnarray}\llabel{As1}
(Y_i -\varepsilon , u^r_s) = \delta_{rs} \delta_{ir} ~~ {\rm for}~~ 
(r,s) \in I, i' \leq i,
\end{eqnarray}
and put
$$
X_{ii} := Y_i - \varepsilon ~~{\rm for}~~ i'\leq i, i=1,\dots,N.
$$
Further, we suppose that
\begin{eqnarray}\llabel{As2}
Z_j, Z_{j'} \in {\rm Lin}~~ \{ Y_i; i' \leq i \} ~~{\rm for}~~  j' < j.
\end{eqnarray}
Then the vector space $\T = {\rm Lin}\, \{ X_{rs}; 
(s,r) \in I \} $ is the quantum tangent space of a left-covariant FODC $\Gamma$ 
on $\cO (O_q(N))$. From the construction and the explicit
 form of the matrix $R$ it is straightforward to check that
that $(X_{rs}, u^i_j) \neq 0 $ if and only if $(r,s)=(j,i)$ for arbitrary indices 
$(s,r) \in I$ and $(i,j) \in I$. This implies that the FODC $\Gamma$ has the 
dimension $N(N+1)/2$ and that the elements $\theta_{ij}= \omega 
(u^i_j), (i,j) \in I,$ form a basis of the space of 
left-invariant one-forms $_{\rm inv}\Gamma$. These facts are in accordance with the 
ordinary first order calculus on the Lie group $O(N)$. Note that the FODC 
$\Gamma$ just constructed depends on the group-like elements $Z_j, Z_{j'}, 
j' < j,$ and $Y_i, i' \le i,$ of $\cO (O_q(N))^\circ$ which can be freely 
chosen such that they satisfy the assumptions (\ref{As1}) and (\ref{As2}). These conditions can be 
easily fulfilled by taking monomials in the main diagonal 
L-functionals $l^{\pm i}_{~~i}$. We make all that more explicit by an example.
\sn 

\noindent {\bf Example: $\cO (O_q(5))$}

Then the $15$ generators of the quantum tangent space $\T$ are 
\begin{eqnarray*}
X_{15} = l^{+1}_{~~5} l^{-1}_{~~1}Z_1,\quad 
X_{25} &=& l^{+1}_{~~4} l^{-1}_{~~1}Z_1,\quad
X_{35} = l^{+1}_{~~3} l^{-1}_{~~1}Z_1,\quad
X_{45} = l^{+1}_{~~2} l^{-1}_{~~1}Z_1,\\
X_{24} &=& l^{+2}_{~~4} l^{-2}_{~~2}Z_2,\quad
X_{34} = l^{+2}_{~~3} l^{-2}_{~~2}Z_2,\\
X_{51} = l^{+5}_{~~1} l^{-5}_{~~5}Z_5,\quad
X_{52} &=& l^{+5}_{~~2} l^{-5}_{~~5}Z_5,\quad
X_{53} = l^{+5}_{~~3} l^{-5}_{~~5}Z_5,\quad
X_{54} = l^{+5}_{~~4} l^{-5}_{~~5}Z_5,\\
X_{42} &=& l^{+4}_{~~2} l^{-4}_{~~4}Z_4,\quad
X_{43} = l^{+4}_{~~3} l^{-4}_{~~4}Z_4,\\
X_{33} = Y_3 - \varepsilon ,\quad
X_{44}&=& Y_4 - \varepsilon ,\quad
X_{55} = Y_5 - \varepsilon 
\end{eqnarray*}
and assumption (\ref{As2}) means that $Z_1, Z_2, Z_4, Z_5 \in 
{\rm Lin}\, \{Y_3, Y_4, Y_5 \}$. \hfill $\Box{}$

\end{document}